\documentclass[a4paper]{amsart}
\usepackage{amscd, amssymb, color, booktabs, chemarrow}
\usepackage{myown}

\usepackage[makeroom]{cancel}
\usepackage{harpoon}

\title[On Theta and Eta Correspondences]{On Theta and Eta Correspondences for Finite Unitary Dual Pairs}
\author{Shu-Yen Pan}
\address{Department of Mathematics,
National Tsing Hua University, Hsinchu 300, Taiwan}
\email{sypan@math.nthu.edu.tw}


\date{\today}

\begin{document}

\begin{abstract}

It is known that the $\Theta$-correspondence for a finite reductive dual pair 
is not one-to-one in general.
In this paper, we propose two maximal one-to-one sub-relations $\underline\theta,\overline\theta$ 
of the $\Theta$-correspondences for a finite dual pair of two unitary groups.
These two correspondences are extensions of the analogue of the $\eta$-correspondence
proposed by Gurevich-Howe.   

\end{abstract}

\maketitle
\tableofcontents


\section{Introduction}

\subsection{}
Let $(\bfG,\bfG')$ be a reductive dual pair over a finite field $\bff_q$ of odd characteristic.
Let $G,G'$ denote the groups of rational points of $\bfG,\bfG'$ respectively.
By restricting the Weil character with respect to a non-trivial character $\psi$ of $\bff_q$
to $G\times G'$, we obtain a decomposition
\[
\omega^\psi_{\bfG,\bfG'}
=\sum_{\rho\in\cale(G),\ \rho'\in\cale(G')}m_{\rho,\rho'}\rho\otimes\rho'
\]
where $m_{\rho,\rho'}=0,1$,
and $\cale(G)$ denotes the set of irreducible characters of $G$.
Define
\begin{align*}
\Theta_{\bfG,\bfG'} &=\{\,(\rho,\rho')\in\cale(G)\times\cale(G')\mid m_{\rho,\rho'}\neq 0\,\};\\
\Theta_{\bfG'}(\rho) &=\{\,\rho'\in\cale(G')\mid (\rho,\rho')\in\Theta_{\bfG,\bfG'}\,\}.
\end{align*}
Now $\Theta_{\bfG,\bfG'}$ establishes a relation between $\cale(G)$ and $\cale(G')$
and will be called the \emph{$\Theta$-correspondence} or \emph{Howe correspondence}.
This correspondence is not one-to-one in general.

When $(\bfG,\bfG')$ is a dual pair of one orthogonal group and one symplectic group and in stable range,
a one-to-one sub-correspondence is defined and called the \emph{$\eta$-correspondence} 
in \cite{gurevich-howe}.
In \cite{pan-eta}, we propose two one-to-one sub-correspondences $\underline\theta,\overline\theta$ 
of $\Theta$ for a general orthogonal/symplectic dual pair $(\bfG,\bfG')$ 
and show that three correspondences $\eta,\underline\theta,\overline\theta$ coincide when
$(\bfG,\bfG')$ is in stable range.
Therefore, both $\underline\theta$ and $\overline\theta$ can be regarded as extensions of 
the $\eta$-correspondence beyond dual pairs in stable range.
The main purpose of this article is to prove the analogous results for dual pairs of two
unitary groups.

\subsection{}
Now let $(\bfG,\bfG')$ be a dual pair of two finite unitary groups.
The one-to-one correspondences $\underline\theta$ and $\overline\theta$ are construct in two steps.
The first step is to define the correspondences on unipotent characters. 
It is well known that the irreducible unipotent characters $\rho_\lambda$ of a unitary group 
$\rmU_n(q)$ are parametrized by the partitions $\lambda$ of $n$.
One the other hand, the unipotent characters are also parametrized by certain symbols $\Lambda_\lambda$
(\cf.~\cite{FS}).
Because the unipotent characters are preserved by the $\Theta$-correspondence,
the unipotent part of the correspondence can be described explicitly in terms of symbols
(\cf.~\cite{amr} and \cite{pan-Lusztig-correspondence}).
Then both correspondences $\underline\theta$ and $\overline\theta$ on unipotent correspondence
are defined in terms of these symbols.

It is known in \cite{amr} that the $\Theta$-correspondence commutes with the 
Lusztig correspondence as follows.
For a semisimple element $s$ is the dual group $G^*$ of $G$,
we can define two subgroups $G^{(1)},G^{(2)}$ so that $C_{G^*}(s)\simeq G^{(1)}\times G^{(2)}$ 
and we have a bijection $\Xi_s\colon\cale(G)_s\rightarrow\cale(G^{(1)}\times G^{(2)})_1$. 
Then the following diagram 
\[
\begin{CD}
\cale(G)_s @> \Theta_{\bfG,\bfG'} >> \cale(G')_{s'} \\
@V \Xi_s VV @VV \Xi_{s'} V \\
\cale(G^{(1)}\times G^{(2)})_1 @> \id\otimes\Theta_{\bfG^{(2)},\bfG'^{(2)}} >> \cale(G'^{(1)}\times G'^{(2)})_1
\end{CD}
\]
is commutative.
Then we can extend both correspondences $\underline\theta$ and $\overline\theta$ outside
the unipotent characters by letting them commute with the Lusztig correspondence.

Then we can show that when $(\bfG,\bfG')=(\rmU_n,\rmU_{n'})$ is in stable range 
(i.e., $n\leq\lfloor\frac{n'}{2}\rfloor$), for $\rho\in\cale(G)$, we have 
\begin{itemize}
\item $\underline\theta(\rho)=\overline\theta(\rho)$;

\item $\overline\theta(\rho)$ is of maximal order in the set $\Theta_{\bfG'}(\rho)$. 
\end{itemize}
This means that both $\underline\theta$ and $\overline\theta$ can be regarded 
as extensions of the analogous $\eta$-correspondence for a dual pair of two unitary groups.

A sub-relation $\vartheta$ of $\Theta$ is called a \emph{theta-relation} of it is 
symmetric, semi-persistent, and compatible with the Lusztig correspondence (\cf.~Subsection~\ref{0402}).
The set of theta-relation can be partially ordered by inclusion.
Then our next result is to show that both $\underline\theta$ and $\overline\theta$ are maximal one-to-one
theta-relation.
This means that any theta-relation which properly contains $\underline\theta$ or $\overline\theta$
will not be one-to-one.
That is, if we require that $\underline\theta,\overline\theta$ to be one-to-one,
then they can not be extended any more.

\subsection{}
The contents of this article is as follows.
In Section 2, we discuss the relation between unipotent characters of a finite unitary group
and certain type of symbols.
In particular, the degree formula of the unipotent character is given in terms of the entries of the symbol.
In section 3, we define the two one-to-one correspondences $\underline\theta$ and $\overline\theta$ on
unipotent characters and show that both correspondences are equal when the dual pair is in stable range.
In the final section, we discuss the relation between the Lusztig correspondence and both correspondences $\underline\theta,\overline\theta$.
Then we show that both $\underline\theta$ and $\overline\theta$ are maximal one-to-one theta-relations.

The strategy of the proofs are similar to those given in \cite{pan-eta} and so some of the proofs
will be sketchy in this article.

\section{Dimension Formula of Unipotent Characters of a Unitary Group}
In this section we give a description of a symbol $\Lambda_\lambda$ associated to a partition $\lambda$.
The we give a formula of the degree of the unipotent character $\rho_\lambda$ in terms of entries of $\Lambda$.
Part of the descriptions in this section is modified from \cite{FS}. 
 
\subsection{Partitions and Unipotent Characters}

Let $\bff_q$ be a finite field of $q$ elements where $q$ is a power of an odd prime $p$.
For a finite set $X$, let $|X|$ denote the cardinality of $X$,
and let $|X|_{p'}$ denote the part of $|X|$ prime to $p$.
It is well-known that
\[
|\rmU_n(q)|_{p'}=\prod_{i=1}^n(q^i-(-1)^i)
\]
where $\rmU_n(q)$ is a finite unitary group defined over $\bff_q$.

Let $\lambda=[\lambda_1,\ldots,\lambda_m]$ denote a \emph{partition} of $n$,
i.e., each $\lambda_i\in\bbN\cup\{0\}$,
$\lambda_1\geq\lambda_2\geq\cdots\geq\lambda_m\geq 0$,
$|\lambda|=\lambda_1+\cdots+\lambda_m=n$ (called the \emph{weight} of $\lambda$).
Each $\lambda_i$ is called a \emph{part} of $\lambda$.
The number $\ell(\lambda)=|\{\,i\mid \lambda_i>0\,\}|$ of positive parts is called the
\emph{length} of $\lambda$.
Two partitions are regarded to be the same if one is obtained from the other by adding several $0$'s,
for example, $[\lambda_1,\ldots,\lambda_m]=[\lambda_1,\ldots,\lambda_m,0]$.
Let $\calp(n)$ denote the set of all partitions of $n$.
For two partitions $\lambda,\lambda'$,
let $\lambda\cup\lambda'$ denote the partitions whose parts consists of the disjoint union of
parts of $\lambda$ and $\lambda'$.

For each index $j$, we define $\lambda^*_j=|\{\,i\mid\lambda_i\geq j\,\}|$.
And we call the partition $\lambda^\rmT=[\lambda_1^*,\lambda_2^*,\ldots]$ the \emph{dual partition} of $\lambda$.

The (Young) \emph{diagram} of $\lambda$ is defined to be the set
of points $(i,j)\in\bbN\times\bbN$ such that $1\leq j\leq\lambda_i$.
To each point $(i,j)$ in the diagram of $\lambda$ we can associate a subset
\[
\xi=\xi_{i,j}=\{\,(i,j')\mid j\leq j'\leq\lambda_i\,\}\cup\{\,(i',j)\mid i\leq i'\leq\lambda^*_i \,\},
\]
which is called the \emph{hook} at $(i,j)$.
The cardinality
\[
|\xi|=\lambda_i+\lambda^*_j-i-j+1.
\]
is called the \emph{hook-length} of $\xi$.
A hook of length $2$ is also called a \emph{$2$-hook}.
A hook $\xi_{i',j'}$ of $\lambda$ is said to \emph{be above} $\xi_{i,j}$ if $j'=j$ and $i'<i$;
a hook $\xi_{i',j'}$ of $\lambda$ is said to \emph{be left to} $\xi_{i,j}$ if $i'=i$ and $j'<j$.

If $\lambda=[\lambda_1,\ldots,\lambda_k]$ ($\lambda_1\geq\cdots\lambda_k\geq 0$) and
$\lambda'=[\lambda'_1,\ldots,\lambda'_k]$ ($\lambda'_1\geq\cdots\lambda'_k\geq 0$) be two partitions
such that either
\begin{itemize}
\item there exists an index $i_0$ such that $\lambda'_{i_0}=\lambda_{i_0}+2$, and
$\lambda'_i=\lambda_i$ for $i\neq i_0$; or

\item there exists an index $i_0$ such that $\lambda'_{i_0}=\lambda_{i_0}+1$, $\lambda'_{i_0+1}=\lambda_{i_0+1}+1$,
and $\lambda'_i=\lambda_i$ for $i\neq i_0,i_0+1$,
\end{itemize}
then we said that \emph{$\lambda'$ is obtained from $\lambda$ by adding a $2$-hook}
or \emph{$\lambda$ is obtained from $\lambda'$ by removing a $2$-hook}.

Let $\lambda$ be a partition of $n$.
After removing all possible $2$-hooks step by step,
the resulting partition is denoted by $\lambda_\infty$ and called the \emph{2-core} of $\lambda$.
Note that $\lambda_\infty=[d,d-1,\ldots,1]$ for some non-negative integer $d$.
A partition $\lambda$ is called \emph{cuspidal} if $\lambda=\lambda_\infty$.

The irreducible unipotent of $\rmU_n(q)$ associated to $\lambda\in\calp(n)$ is denoted by $\rho_\lambda$.
It is known that the dimensional of $\rho_\lambda$ is  
\[
\dim(\rho_\lambda)=q^{\kappa(\lambda)}g_\lambda(q)
\]
where $\kappa(\lambda)=\sum_{i=1}^{\ell(\lambda)}(i-1)\lambda_i$ and 
$g_\lambda$ is a polynomial given by
\[
g_\lambda
=g_\lambda(q)
=\frac{\prod_{i=1}^n(q^i-(-1)^i)}{\prod_{\xi}(q^{|\xi|}-(-1)^{|\xi|})}
\]
where $\prod_{\xi}$ runs over all hooks $\xi$ of $\lambda$.
It is clear that $g_\lambda=g_{\lambda^\rmT}$.
Now $\dim(\rho_\lambda)$ is a polynomial of $q$ whose degree is denoted by 
$\ord(\rho_\lambda)$.
From above, we know that 
\[
\ord(\rho_\lambda)=\kappa(\lambda)+\frac{n(n+1)}{2}-\sum_\xi|\xi|.
\]

\begin{exam}\label{0202}
Suppose that $\lambda=[d,d-1,\ldots,1]$.
Then $\ell(\lambda)=d$ and
\[
\kappa(\lambda)
=(d-1)+2(d-2)+\cdots+(d-1)\cdot1
=\binom{d}{2}+\binom{d-1}{2}+\cdots+\binom{2}{2}.
\]
Now $\lambda$ has exactly $i$ hooks of length $2d-(2i-1)$ for $i=1,2,\ldots,d$,
and no hooks of other lengths,
so
\[
\dim(\rho_\lambda)
=\frac{q^{\binom{d}{2}+\binom{d-1}{2}+\cdots+\binom{2}{2}}\prod_{i=1}^n(q^i-(-1)^i)}
{(q+1)^{d}(q^3+1)^{d-1}\cdots(q^{2d-3}+1)^2(q^{2d-1}+1)}.
\]
It is known that $\rho_\lambda$ is cuspidal and the dimension of $\rho_\lambda$ is also given in
\cite{lg}.
\end{exam}

\subsection{Partitions and $\beta$-Sets}
A \emph{$\beta$-set} $X=\{x_1,\ldots,x_m\}$ is a finite set of non-negative integers with elements
written in decreasing order, i.e., $x_1>x_2>\cdots>x_m$.
Define an equivalence relation (denoted by ``$\sim$'') on $\beta$-sets generated by
\[
X\sim \{\,x+1\mid x\in X\}\cup\{0\}.
\]
Define the \emph{rank} of $X$ by
\[
{\rm rk}(X)=\sum_{x\in X} x-\binom{|X|}{2}.
\]
It is clear that ${\rm rk}(X')={\rm rk}(X)$ if $X'\sim X$
and the mapping
\begin{equation}\label{0213}
\Upsilon\colon\{x_1,x_2,\ldots,x_m\}\mapsto[x_1-(m-1),x_2-(m-2),\ldots,x_{m-1}-1,x_m]
\end{equation}
gives a bijection from set of equivalence classes of $\beta$-sets of rank $n$ onto the set $\calp(n)$
of partitions of $n$.

For a $\beta$-set $X$, we define four $\beta$-sets:
\begin{align*}
X^1 &=\{\,x\in X\mid x\equiv 1\pmod 2\,\}, & \overline X^1 &=\left\{\,\frac{x-1}{2}\mid x\in X^1\,\right\};\\
X^0 &=\{\,x\in X\mid x\equiv 0\pmod 2\,\}, & \overline X^0 &=\left\{\,\frac{x}{2}\mid x\in X^0\,\right\} .
\end{align*}

\begin{lem}
Let $X,X'$ be two $\beta$-sets.
If $X'\sim X$, 
then the difference $|X'^1|-|X'^0|$ is equal to either $|X^1|-|X^0|$ or $|X^0|-|X^1|-1$.
\end{lem}
\begin{proof}
Suppose that $X'=\{\,x+1\mid x\in X\,\}\cup\{0\}$.
Then clearly $|X'^1|=|X^0|$ and $|X'^0|=|X^1|+1$, and hence
$|X'^1|-|X'^0|=|X^0|-|X^1|-1$.
Note that $|X'^0|-|X'^1|-1=|X^1|-|X^0|$ and the lemma is proved.
\end{proof}

\begin{lem}
Let $X$ be a $\beta$-set of rank $n$, and let $\delta=|X^1|-|X^0|$.
Then $n-\frac{\delta(\delta+1)}{2}$ is a non-negative even integer.
\end{lem}
\begin{proof}
Suppose that $|X^1|=m_1$ and $|X^0|=m_0$.
Then $|X|=m_0+m_1$ and $\delta=m_1-m_0$.
Now
\begin{align*}
& n-\frac{(m_1-m_0)(m_1-m_0+1)}{2} \\
&= \sum_{x\in X}x-\frac{(m_1+m_0)(m_1+m_0-1)}{2}-\frac{(m_1-m_0)(m_1-m_0+1)}{2} \\
&= \sum_{x\in X}x-m_1^2-m_0(m_0+1)
\end{align*}
Because any element in $X^0$ is even and $m_0(m_0+1)$ is also even,
we have
\[
\sum_{x\in X}x\equiv \sum_{x\in X^1}x \equiv m_1 \equiv m_1^2\pmod 2.
\]
So the lemma is proved.
\end{proof}

For a partition $\lambda$, we define a $\beta$-set
\begin{align*}
X=X_\lambda &=\begin{cases}
\{\lambda_1+(m-1),\lambda_2+(m-2),\ldots,\lambda_{m-1}+1,\lambda_m\}, & \text{if $\ell(\lambda)+\ell(\lambda_\infty)$ is even};\\
\{\lambda_1+m,\lambda_2+(m-1),\ldots,\lambda_{m-1}+2,\lambda_m+1,0\}, & \text{otherwise},
\end{cases}
\end{align*}
It is clear that the rank of $X_\lambda$ is equal to $|\lambda|$.
Moreover, $X_{\lambda'}\sim X_\lambda$ if $\lambda'$ is obtained from $\lambda$
by adding several $0$'s.

\begin{lem}
Let $\lambda$ be a partition.
Then $|X_\lambda^1|\geq |X_\lambda^0|$.
\end{lem}
\begin{proof}
Let $X=X_\lambda$.
After removing all possible $2$-hooks,
we will obtain a $\beta$-set $X_\infty$ of the form
\[
\{2k-1,2k-3,\ldots,1\}\cup\{2l-2,2l-4,\ldots,0\}
\]
for some non-negative integers $k,l$, and $\lambda_\infty=[d,d-1,\ldots,1]$ where
\[
d=\begin{cases}
k-l, & \text{if $k\geq l$};\\
l-k-1, & \text{if $k<l$}.
\end{cases}
\]
Now $|X^1|=k$, $|X^0|=l$ and $\ell(\lambda_\infty)=d$.

\begin{enumerate}
\item Suppose that $\ell(\lambda)+\ell(\lambda_\infty)$ is even.
Then $\ell(\lambda)=|X|=k+l$.
If $l>k$, then $\ell(\lambda_\infty)+\ell(\lambda)=(l-k-1)+(l+k)=2l-1$, which is odd,
and we get a contradiction.

\item Suppose that $\ell(\lambda)+\ell(\lambda_\infty)$ is odd.
Then $\ell(\lambda)=|X|-1=l+k-1$.
If $l>k$, then $\ell(\lambda_\infty)+\ell(\lambda)=(l-k-1)+(l+k-1)=2l-2$, which is even,
and we get a contradiction, again.
\end{enumerate}
Therefore, we must have $|X^1|\geq |X^0|$.
\end{proof}

\begin{lem}\label{0211}
Let $\lambda$ be a partition.
Then $\ell(\lambda_\infty)=|X_\lambda^1|-|X^0_\lambda|$.
\end{lem}
\begin{proof}
Let $X=X_\lambda$.
After removing all possible $2$-hooks,
we will obtain a $\beta$-set $X_\infty$ of the form
\[
\{2k-1,2k-3,\ldots,1\}\cup\{2l-2,2l-4,\ldots,0\}
\]
for some non-negative integers $k,l$.
From the previous lemma, we know that $k\geq l$ and $\ell(\lambda_\infty)=k-l=|X^1|-|X^0|$.
\end{proof}

For a partition $\lambda$ (or for a $\beta$-set $X_\lambda$),
we associate a symbol
\[
\Lambda_\lambda=
\begin{cases}
\binom{X^1_\lambda}{X^0_\lambda}, & \text{if $\ell(\lambda_\infty)$ is even};\\
\binom{X^0_\lambda}{X^1_\lambda}, & \text{if $\ell(\lambda_\infty)$ is odd}.
\end{cases}
\]
For a symbol $\Lambda$, let $\Lambda^*$ (resp.~$\Lambda_*$) denote the first row (resp.~second row)
of $\Lambda$. 
Recall that the \emph{defect} of a symbol $\Lambda$ is defined to be
${\rm def}(\Lambda)=|\Lambda^*|-|\Lambda_*|$.
Then by Lemma~\ref{0211}, we have
\begin{equation}\label{0214}
{\rm def}(\Lambda_\lambda)=\begin{cases}
\phantom{-}\ell(\lambda_\infty), & \text{if $\ell(\lambda_\infty)$ is even};\\
-\ell(\lambda_\infty), & \text{if $\ell(\lambda_\infty)$ is odd}.
\end{cases}
\end{equation}

For a non-negative even integer or a negative odd integer $\delta$,
let $\cals_{n,\delta}$ denote the set of equivalence class of symbols $\Lambda_\lambda$
such that ${\rm rk}(X_\lambda)=n$ and ${\rm def}(\Lambda_\lambda)=\delta$.
Let $\calp_2(n)$ denote the set of \emph{bi-partitions} of $n$, i.e.,
the set of $\sqbinom{\mu}{\nu}$ where $\mu,\nu$ are partitions such that $|\sqbinom{\mu}{\nu}|=|\mu|+|\nu|=n$.
Now we define a mapping
\begin{equation}\label{0212}
\Lambda_\lambda\mapsto\begin{cases}
\sqbinom{\Upsilon(\overline X^1)}{\Upsilon(\overline X^0)}, & \text{if $\ell(\lambda_\infty)$ is even};\\
\sqbinom{\Upsilon(\overline X^0)}{\Upsilon(\overline X^1)}, & \text{if $\ell(\lambda_\infty)$ is odd}
\end{cases}
\end{equation}
where $\Upsilon$ is given in (\ref{0213}).
By abusing the notation a little bit,
the above mapping is also denoted by $\Upsilon$.
It is easy to check that (\cf.~\cite{FS} p.223)
\begin{equation}\label{0516}
|\lambda|=|\lambda_\infty|+2|\Upsilon(\Lambda_\lambda)|.
\end{equation}
Moreover, $\Upsilon$ gives a bijection from $\cals_{n,\delta}$ onto
$\calp_2(\frac{1}{2}(n-\frac{|\delta|(|\delta|+1)}{2}))$.

\begin{exam}
Suppose that $\lambda=[d,d-1,\ldots,1]$.
Then $\lambda_\infty=\lambda$ and hence
$X_\lambda=\{2d-1,2d-3,\ldots,1\}$.
Then $\overline X^1=\{d-1,d-2,\ldots,0\}$, $\overline X^0=\emptyset$, and hence
\[
\Upsilon(\Lambda_\lambda)=\sqbinom{0}{0}\in\calp_2(0).
\]
\end{exam}

\begin{exam}
Suppose that $\lambda=[1,\ldots,1]\in\calp(n)$.
\begin{enumerate}
\item Suppose that $n$ is even.
Then $\lambda_\infty=[0]$, and $\ell(\lambda)+\ell(\lambda_\infty)=n$.
Then $X_\lambda=\{n,n-1,\ldots,1\}$,
and hence $\overline X^1=\{\frac{n}{2}-1,\frac{n}{2}-2,\ldots,0\}$ and $\overline X^0=\{\frac{n}{2},\frac{n}{2}-1,\ldots,1\}$.
Then
\[
\Upsilon(\Lambda_\lambda)=\sqbinom{0}{1,1,\ldots,1}\in\calp_2(\tfrac{n}{2}).
\]

\item Suppose that $n$ is odd.
Then $\lambda_\infty=[1]$, and $\ell(\lambda)+\ell(\lambda_\infty)=n+1$.
Then $X_\lambda=\{n,n-1,\ldots,1\}$,
and hence $\overline X^1=\{\frac{n-1}{2},\frac{n-3}{2},\ldots,0\}$ and
$\overline X^0=\{\frac{n-1}{2},\frac{n-3}{2},\ldots,1\}$.
Then
\[
\Upsilon(\Lambda_\lambda)=\sqbinom{1,1,\ldots,1}{0}\in\calp_2(\tfrac{n-1}{2}).
\]
\end{enumerate}
\end{exam}

Define
\[
\cals_{\rmU_n}=\{\,\Lambda_\lambda\mid\lambda\in\calp(n)\,\}.
\]
Therefore, $\cals_{n,\delta}\subset\cals_{\rmU_n}$ if and only if
\begin{itemize}
\item $\delta$ is either a non-negative integer or a negative odd integer; and

\item $\frac{1}{2}(n-\frac{|\delta|(|\delta|+1)}{2})$ is a non-negative integer.
\end{itemize}

\begin{exam}
From above, we see that
\[
\cals_{\rmU_7}=\cals_{7,-1}\cup\cals_{7,2},\qquad
\cals_{\rmU_8}=\cals_{8,0}\cup\cals_{8,-3},\qquad
\cals_{\rmU_{10}}=\cals_{10,0}\cup\cals_{10,-3}\cup\cals_{10,4}.
\]
\end{exam}

\begin{rem}
The definition of the symbol $\Lambda_\lambda\in\cals_{\rmU_n}$ associated to a given partition $\lambda$ of $n$ 
is slightly different from that given in \cite{pan-Lusztig-correspondence}.
The new definition here is more convenient for us, 
in particular, (\ref{0214}) is simpler than lemma 5.5 in \cite{pan-Lusztig-correspondence}.
\end{rem}

\subsection{Order of a $\beta$-Set}

For two $\beta$-sets $A,B$, we define several polynomials in $q$:
\begin{align*}
\Delta(A)
&:= \prod_{a,a'\in A,\ a>a'}(q^a-q^{a'}) \\
\Theta(A)
&:= \prod_{a\in A}\prod_{h=1}^a(q^a-(-1)^a) \\
\Xi(A,B)
&:= \prod_{a\in A,\ b\in B}(q^a+q^b) \\
f_{A,B}
&:=\frac{\Delta(A)\Delta(B)\Xi(A,B)}{\Theta(A)\Theta(B)q^{\binom{|A|+|B|-1}{2}+\binom{|A|+|B|-2}{2}+\cdots+\binom{2}{2}}}.
\end{align*}
For a $\beta$-set $X$,
we define $f_X=f_{X^0,X^1}$.

\begin{lem}
If $X\sim X'$,
then $f_X=f_{X'}$.
\end{lem}
\begin{proof}
Without loss of generality, we may assume that
\[
X'=\{\, x+1\mid x\in X\,\}\cup\{0\}.
\]
Then
\[
X'^0=\{\,b+1\mid b\in X^1\,\}\cup\{0\}\quad\text{ and }\quad
X'^1=\{\,a+1\mid a\in X^0\,\}.
\]
From definitions above, it is clear that
\begin{align*}
\Delta(X'^0)
&= \Delta(X^1)\cdot q^{\binom{|X^1|}{2}}\cdot\prod_{b\in X^1}(q^{b+1}-1) \\
\Delta(X'^1)
&= \Delta(X^0)\cdot q^{\binom{|X^0|}{2}} \\
\Theta(X'^0)
= \Theta(X^1)\cdot\prod_{b\in X^1}(q^{b+1}-(-1)^{b+1})
&= \Theta(X^1)\cdot\prod_{b\in X^1}(q^{b+1}-1) \\
\Theta(X'^1)
= \Theta(X^0)\cdot\prod_{a\in X^0}(q^{a+1}-(-1)^{a+1})
&= \Theta(X^0)\cdot\prod_{a\in X^0}(q^{a+1}+1).
\end{align*}
Moreover, we have
\begin{align*}
\Xi(X'^0,X'^1)
&= \prod_{a'\in X'^0,\ b'\in X'^1}(q^{a'}+q^{b'}) \\
&= \prod_{a\in X^0,\ b\in X^1}(q^{a+1}+q^{b+1})\cdot\prod_{a\in X^0}(q^{a+1}+1) \\
&= \Xi(X^0,X^1)\cdot q^{|X^0||X^1|}\cdot\prod_{a\in X^0}(q^{a+1}+1)
\end{align*}
Now $|X'^0|=|X^1|+1$, $|X'^1|=|X^0|$, and
\[
\binom{|X'|-1}{2}=\binom{|X^0|+|X^1|}{2}=\binom{|X^0|}{2}+\binom{|X^1|}{2}+|X^0|\cdot|X^1|.
\]
Then the lemma follows.
\end{proof}

Now we want to express the dimension of $\rho_\lambda$ in terms of the $\beta$-set $X_\lambda$.
Similar but different expression can be found in \cite{FS}.

\begin{prop}
Let $\lambda$ be a partition of $n$, and let $X=X_\lambda$.
Then
\begin{align*}
\dim(\rho_\lambda)
&= \frac{\Delta(X^0)\Delta(X^1)\Xi(X^0,X^1)}
{\Theta(X^0)\Theta(X^1)q^{\binom{|X|-1}{2}+\binom{|X|-2}{2}+\cdots+\binom{2}{2}}}\cdot \prod_{i=1}^n(q^i-(-1)^i) \\
&= f_X\cdot |\rmU_n(q)|_{p'}.
\end{align*}
\end{prop}
\begin{proof}
Note that every partition $\lambda$ is built from a cuspidal partition $\lambda_\infty$
by adding several $2$-hooks.
So we now prove the proposition by induction on the number of $2$-hooks.

First suppose that $\lambda=\lambda_\infty=[d,d-1,\ldots,1]$ for some $d$.
Then $n=\frac{d(d+1)}{2}$,
$X=X^1=\{2d-1,2d-3,\ldots,1\}$ and $X^0=\emptyset$.
Therefore $|X|=|X^1|=d$, and
\begin{align*}
\Delta(X^0)=\Theta(X^0)=\Xi(X^0,X^1)
&=1, \\
\Delta(X^1)
&=q^{\binom{d}{2}+2\left[\binom{d-1}{2}+\binom{d-2}{2}+\cdots+\binom{2}{2}\right]}\prod_{x=1}^{d-1}\prod_{h=1}^{x}(q^{2h}-1), \\
\Theta(X^1)
&= (q+1)^{d}\prod_{x=1}^{d-1}\prod_{h=1}^x \left[(q^{2h}-1)(q^{2h+1}+1)\right].
\end{align*}
So we see that
\[
f_X\cdot|\rmU_n(q)|_{p'}
=\frac{q^{\binom{d}{2}+\binom{d-1}{2}+\cdots+\binom{2}{2}}\prod_{i=1}^n(q^i-(-1)^i)}
{(q+1)^{d}(q^3+1)^{d-1}\cdots(q^{2d-3}+1)^2(q^{2d-1}+1)}
=\dim(\rho_\lambda)
\]
from Example~\ref{0202}.
Therefore the lemma is true when the partition $\lambda$ is cuspidal.

Next we assume that the lemma is true for some partition $\lambda$ of $n$ and suppose that
$\lambda'$ is obtained from $\lambda$ by adding a $2$-hook.
Then $\lambda'$ is a partition of $n+2$.
We have the following two situations:
\begin{enumerate}
\item There is an index $i_0$ such that $\lambda'_{i_0}=\lambda_{i_0}+2$
and $\lambda_i'=\lambda_i$ for any $i\neq i_0$.
Moreover, without loss of generality, we may assume that the
part of rows below $\lambda_{i_0}$ is cuspidal, i.e.,
\begin{equation}\label{0201}
[\lambda_{i_0+1},\lambda_{i_0+2},\ldots,\lambda_{\ell(\lambda)}]=[d',d'-1,\ldots,1]
\end{equation}
for some non-negative integer $d'$.
Then clearly $\kappa(\lambda')=\kappa(\lambda)+2(i_0-1)$,
and
\begin{multline*}
\frac{g_{\lambda'}}{g_{\lambda}}
=\frac{(q^{n+1}-(-1)^{n+1})(q^{n+2}-(-1)^{n+2})}{(q-1)(q^2+1)}\cdot
\prod_{\xi'}\frac{q^{|\xi'|}-(-1)^{|\xi'|}}{q^{|\xi'|+2}-(-1)^{|\xi'|}} \\
\cdot \prod_{\xi''}\frac{q^{|\xi''|}-(-1)^{|\xi''|}}{q^{|\xi''|+1}-(-1)^{|\xi''|+1}}\cdot
\prod_{\xi'''}\frac{q^{|\xi'''|}-(-1)^{|\xi'''|}}{q^{|\xi'''|+1}-(-1)^{|\xi'''|+1}}
\end{multline*}
where $\prod_{\xi'}$ runs over all hooks left to $(i_0,\lambda_{i_0})$,
$\prod_{\xi''}$ runs over all hooks above $(i_0,\lambda_{i_0}+1)$,
and $\prod_{\xi'''}$ runs over all hooks above $(i_0,\lambda_{i_0}+2)$.
Note that
\[
|\xi_{i,\lambda_{i_0+1}}|=|\xi_{i,\lambda_{i_0+2}}|+1
\]
for $i=1,\ldots,i_0-1$, so
\begin{multline}\label{0203}
\frac{g_{\lambda'}}{g_{\lambda}}
=\frac{(q^{n+1}-(-1)^{n+1})(q^{n+2}-(-1)^{n+2})}{(q-1)(q^2+1)}\cdot
\prod_{\xi'}\frac{q^{|\xi'|}-(-1)^{|\xi'|}}{q^{|\xi'|+2}-(-1)^{|\xi'|}} \\
\cdot \prod_{\xi'''}\frac{q^{|\xi'''|}-(-1)^{|\xi'''|}}{q^{|\xi'''|+2}-(-1)^{|\xi'''|}}.
\end{multline}
Now we have two possibilities:
\begin{enumerate}
\item $X^1_{\lambda'}=X^1_\lambda$ and there is an $a_0\in X^0_\lambda$ such that
$a_0+2\not\in X^0_\lambda$ and $X^0_{\lambda'}=X^0_\lambda\cup\{a_0+2\}\smallsetminus\{a_0\}$.
Moreover, we know that if $b\in X^1_\lambda$, then either $b>a_0+2$ or $b<a_0$.
Then we have
\begin{align}\label{0204}
\begin{split}
\Delta(X^1_{\lambda'}) &= \Delta(X^1_\lambda); \\
\Theta(X^1_{\lambda'}) &= \Theta(X^1_\lambda); \\
\Theta(X^0_{\lambda'}) &= \Theta(X^0_\lambda)(q^{a_0+1}+1)(q^{a_0+2}-1);
\end{split}
\end{align}
\begin{align}\label{0205}
\begin{split}
\frac{\Delta(X^0_{\lambda'})}{\Delta(X^0_\lambda)}
&= \prod_{a\in X^0_\lambda,\ a>a_0+2}\frac{q^a-q^{a_0+2}}{q^a-q^{a_0}}
\cdot \prod_{a\in X^0_\lambda,\ a<a_0}\frac{q^{a_0+2}-q^a}{q^{a_0}-q^a} \\
&= q^{2l}
\cdot \prod_{a\in X^0_\lambda,\ a>a_0+2}\frac{q^{a-a_0-2}-1}{q^{a-a_0}-1}
\cdot \prod_{a\in X^0_\lambda,\ a<a_0}\frac{q^{a_0+2-a}-1}{q^{a_0-a}-1}
\end{split}
\end{align}
where $l$ is the number of elements $a\in X^0_\lambda$ which is greater than $a_0$;
\begin{align}\label{0206}
\begin{split}
\frac{\Xi(X^0_{\lambda'},X^1_{\lambda'})}{\Xi(X^0_\lambda,X^1_\lambda)}
&= \prod_{b\in X^1_\lambda}\frac{q^{a_0+2}+q^b}{q^{a_0}+q^b} \\
&= q^{2l'}
\cdot \prod_{b\in X^1_\lambda,\ b>a_0+2}\frac{q^{b-a_0-2}+1}{q^{b-a_0}+1}
\cdot \prod_{b\in X^1_\lambda,\ b<a_0}\frac{q^{a_0-b+2}+1}{q^{a_0-b}+1}
\end{split}
\end{align}
where $l'$ is the number of elements $b\in X^1_\lambda$ which is greater than $a_0$.
Now $l+l'$, the number of rows $\lambda_i$ which is greater than $\lambda_{i_0}$,
i.e., $l+l'=i_0-1$.
Now (\ref{0201}) implies that
\[
\{\,a\in X^0_\lambda\mid a<a_0\,\}=\emptyset\quad\text{ and }\quad
\{\,b\in X^1_\lambda\mid b<a_0\,\}=\{2d'-1,2d'-3,\ldots,1\}.
\]
Then
\[
\prod_{a\in X^0_\lambda,\ a<a_0}\frac{q^{a_0+2-a}-1}{q^{a_0-a}-1}
\cdot\prod_{b\in X^1_\lambda,\ b<a_0}\frac{q^{a_0-b+2}+1}{q^{a_0-b}+1}
=\frac{q^{a_0+1}+1}{q^{a_0-2d'+1}+1}.
\]

There is a one-to-one correspondence between the set of
hooks $\xi'''$ above $(i_0,\lambda_{i_0}+2)$ and the set
\[
\{\,a\in X^0_\lambda\mid a>a_0+2\,\}\cup\{\,b\in X^1_\lambda\mid b>a_0+2\,\}
\]
and
\[
|\xi'''|=a-a_0-2\quad\text{ or }\quad
|\xi'''|=b-a_0-2.
\]
Note that $a-a_0-2$ is even and $b-a_0-2$ is odd.
Therefore
\begin{equation}\label{0207}
\prod_{\xi'''}\frac{q^{|\xi'''|}-(-1)^{|\xi'''|}}{q^{|\xi'''|+2}-(-1)^{|\xi'''|}}
=\prod_{a\in X^0_\lambda,\ a>a_0+2}\frac{q^{a-a_0-2}-1}{q^{a-a_0}-1}
\cdot\prod_{b\in X^1_\lambda,\ b>a_0+2}\frac{q^{b-a_0-2}+1}{q^{b-a_0}+1}.
\end{equation}
Now consider the hooks left to $(i_0,\lambda_{i_0})$.
From (\ref{0201}), we know that
\begin{equation}\label{0208}
\prod_{\xi'}\frac{q^{|\xi'|}-(-1)^{|\xi'|}}{q^{|\xi'|+2}-(-1)^{|\xi'|}}
=\frac{(q-1)(q^2+1)}{(q^{a_0+2}-1)(q^{a_0-2d'+1}+1)}.
\end{equation}
Therefore, by (\ref{0203}), (\ref{0204}), (\ref{0205}), (\ref{0206}), (\ref{0207}), and (\ref{0208}),
we conclude that
\[
\frac{f_{X^0_{\lambda'},X^1_{\lambda'}}|\rmU_{n+2}(q)|_{p'}}{f_{X^0_\lambda,X^1_\lambda}|\rmU_n(q)|_{p'}}
=\frac{q^{\kappa(\lambda')}g_{\lambda'}}{q^{\kappa(\lambda)}g_\lambda}
=\frac{\deg(\rho_{\lambda'})}{\deg(\rho_\lambda)}.
\]
Thus the proposition is true for $\lambda'$ by induction hypothesis.

\item $X^0_{\lambda'}=X^0_\lambda$ and there is a $b_0\in X^1_\lambda$ such that $b_0+2\not\in X^1_\lambda$
and $X^1_{\lambda'}=X^1_\lambda\cup\{b_0+2\}\smallsetminus\{b_0\}$.
The proof is similar to case (a) and is skipped.
\end{enumerate}

\item There exists an index $i_0$ such that $\lambda_{i_0}'=\lambda_{i_0}+1$ and
$\lambda'_{i_0+1}=\lambda_{i_0+1}$,
and $\lambda_j'=\lambda_j$ for any $j\neq i_0,i_0+1$.
Again, without loss of generality, we may assume that the
part of rows below $\lambda_{i_0+1}$ is cuspidal, i.e.,
\begin{equation}
[\lambda_{i_0+2},\lambda_{i_0+3},\ldots,\lambda_{\ell(\lambda)}]=[d',d'-1,\ldots,1]
\end{equation}
for some non-negative integer $d'$.
Then $\kappa(\lambda')=\kappa(\lambda)+i_0+(i_0-1)=\kappa(\lambda)+(2i_0-1)$, and
\begin{multline*}
\frac{g_{\lambda'}}{g_{\lambda}}
=\frac{(q^{n+1}-(-1)^{n+1})(q^{n+2}-(-1)^{n+2})}{(q-1)(q^2+1)}\cdot
\prod_{\xi'}\frac{q^{|\xi'|}-(-1)^{|\xi'|}}{q^{|\xi'|+1}-(-1)^{|\xi'|+1}} \\
\cdot \prod_{\xi''}\frac{q^{|\xi''|}-(-1)^{|\xi''|}}{q^{|\xi''|+1}-(-1)^{|\xi''|+1}}\cdot
\prod_{\xi'''}\frac{q^{|\xi'''|}-(-1)^{|\xi'''|}}{q^{|\xi'''|+2}-(-1)^{|\xi'''|}}
\end{multline*}
where $\prod_{\xi'}$ runs over all hooks left to $(i_0,\lambda_{i_0})$,
$\prod_{\xi''}$ runs over all hooks left to $(i_0+1,\lambda_{i_0+1})$,
and $\prod_{\xi'''}$ runs over all hooks above $(i_0,\lambda_{i_0}+1)$.
Note that
\[
|\xi_{i_0,j}|=|\xi_{i_0+1,j}|+1
\]
for $j=1,\ldots,\lambda_{i_0}$, so
\begin{multline}
\frac{g_{\lambda'}}{g_{\lambda}}
=\frac{(q^{n+1}-(-1)^{n+1})(q^{n+2}-(-1)^{n+2})}{(q-1)(q^2+1)}\cdot
\prod_{\xi''}\frac{q^{|\xi''|}-(-1)^{|\xi''|}}{q^{|\xi''|+2}-(-1)^{|\xi''|}} \\
\cdot \prod_{\xi'''}\frac{q^{|\xi'''|}-(-1)^{|\xi'''|}}{q^{|\xi'''|+2}-(-1)^{|\xi'''|}}.
\end{multline}

\begin{enumerate}
\item $X^1_{\lambda'}=X^1_\lambda$ and there is an $a_0\in X^0_\lambda$ such that
$a_0+2\not\in X^0_\lambda$ and $X^0_{\lambda'}=X^0_\lambda\cup\{a_0+2\}\smallsetminus\{a_0\}$.
Moreover, we know that there exists a (unique) element $b_0\in X^1_\lambda$ such that $b_0=a_0+1$.
Then we have
\begin{align*}
\Delta(X^1_{\lambda'}) &= \Delta(X^1_\lambda); \\
\Theta(X^1_{\lambda'}) &= \Theta(X^1_\lambda); \\
\Theta(X^0_{\lambda'}) &= \Theta(X^0_\lambda)(q^{a_0+1}+1)(q^{a_0+2}-1);
\end{align*}
\begin{align*}
\frac{\Delta(X^0_{\lambda'})}{\Delta(X^0_\lambda)}
&= q^{2l} \cdot \prod_{a\in X^0_\lambda,\ a>a_0+2}\frac{q^{a-a_0-2}-1}{q^{a-a_0}-1}
\cdot \prod_{a\in X^0_\lambda,\ a<a_0}\frac{q^{a_0+2-a}-1}{q^{a_0-a}-1}
\end{align*}
where $l$ is the number of elements $a\in X^0_\lambda$ which is greater than $a_0$.
Note that there is an element $b_0\in X^1_\lambda$ such that $b_0=a_0+1$,
so
\[
\frac{q^{a_0+2}+q^{b_0}}{q^{a_0}+q^{b_0}}=q.
\]
Then
\begin{align*}
\frac{\Xi(X^0_{\lambda'},X^1_{\lambda'})}{\Xi(X^0_\lambda,X^1_\lambda)}
&= q^{2l'-1}\cdot \prod_{b\in X^1_\lambda,\ b>a_0+2}\frac{q^{b-a_0-2}+1}{q^{b-a_0}+1}
\cdot \prod_{b\in X^1_\lambda,\ b<a_0}\frac{q^{a_0-b+2}+1}{q^{a_0-b}+1}
\end{align*}
where $l'$ is the number of elements $b\in X^1_\lambda$ which is greater than $a_0$.
Now $l+l'$, the number of rows $\lambda_i$ which is greater than $\lambda_{i_0+1}$,
i.e., $l+l'=i_0$.
Then $\kappa(\lambda')-\kappa(\lambda)=2i_0-1=2(l+l')-1$.
Again, we have
\begin{align*}
\prod_{\xi'''}\frac{q^{|\xi'''|}-(-1)^{|\xi'''|}}{q^{|\xi'''|+2}-(-1)^{|\xi'''|}}
&= \prod_{a\in X^0_\lambda,\ a>a_0+2}\frac{q^{a-a_0-2}-1}{q^{a-a_0}-1}
\cdot\prod_{b\in X^1_\lambda,\ b>a_0+2}\frac{q^{b-a_0-2}+1}{q^{b-a_0}+1} \\
\prod_{\xi''}\frac{q^{|\xi''|}-(-1)^{|\xi''|}}{q^{|\xi''|+2}-(-1)^{|\xi''|}}
&= \frac{(q-1)(q^2+1)}{(q^{a_0+2}-1)(q^{a_0-2d'+1}+1)}.
\end{align*}
Again, the proposition is true for $\lambda'$ by induction hypothesis.

\item $X^0_{\lambda'}=X^0_\lambda$ and there is a $b_0\in X^1_\lambda$ such that $b_0+2\not\in X^1_\lambda$
and $X^1_{\lambda'}=X^1_\lambda\cup\{b_0+2\}\smallsetminus\{b_0\}$.
The proof is similar to case (a) and is skipped.
\end{enumerate}
\end{enumerate}
Finally the proposition is proved for all cases.
\end{proof}

For $\lambda\in\calp(n)$, 
let $\ord(X_\lambda)$ or $\ord(\Lambda_\lambda)$ also denote the degree of the polynomial $\dim(\rho_\lambda)$ in $q$.

\begin{lem}\label{0209}
Let $\lambda\in\calp(n)$.
Suppose that $X_\lambda=\{x_1,\ldots,x_m\}$ where $x_1>x_2>\cdots>x_m$.
Then
\[
\ord(X_\lambda)=\sum_{i=1}^m(m-i)x_i-\sum_{i=1}^m x_i(x_i+1)+n(n+1)-\frac{m(m-1)(m-2)}{6}.
\]
\end{lem}
\begin{proof}
Clearly, we have
\begin{align*}
\deg(\Delta(X^0)\Delta(X^1)\Xi(X^0,X^1))
&= \sum_{i=1}^{m}(m-i)x_i \\
\deg(\Theta(X^0)\Theta(X^1))
&= \sum_{i=1}^{m}x_i(x_i+1) \\
\deg(\binom{m-1}{2}+\cdots+\binom{2}{2})
&= \frac{m(m-1)(m-2)}{6} \\
\deg(|\rmU_n(q)|_{p'})
&= n(n+1).
\end{align*}
Then the lemma follows.
\end{proof}

\begin{lem}\label{0210}
Let $X=\{x_1,\ldots,x_m\}$ and $X'=\{x'_1,\ldots,x'_m\}$ be two $\beta$-sets.
Suppose that there are two indices $k<l$ such that $x'_k=x_k+1$, $x'_l=x_l-1$,
and $x'_i=x_i$ for $i\neq k,l$.
Then $\ord(X)>\ord(X')$.
\end{lem}
\begin{proof}
We have
\begin{align*}
(m-k)(x_k+1)-(x_k+1)(x_k+2)-(m-k)x_k+x_k(x_k+1)
&=m-k-2x_k-2, \\
(m-l)(x_l-1)-(x_l-1)x_l-(m-l)x_l+x_l(x_l+1)
&=l-m+2x_l.
\end{align*}
Because $x_1>x_2>\cdots>x_m$, we must have $x_k-x_l\geq k-l$,
and hence by Lemma~\ref{0209}
\[
\ord(X')-\ord(X)=l-k+2(x_l-x_k)-2<0.
\]
\end{proof}

\section{Two Correspondences on Unipotent Characters}
In this section,
we consider a dual pair of two unitary groups, i.e.,
$(\bfG,\bfG')=(\rmU_n,\rmU_{n'})$ for some non-negative $n,n'$.

\subsection{Theta correspondence}
Let $\lambda=[\lambda_1,\ldots,\lambda_m]$ and $\lambda'=[\lambda_1,\ldots,\lambda'_{m'}]$
be two partitions.
By adding some $0$'s if necessary,
we may assume that $m=m'$.
We say that
\[
\lambda\preccurlyeq\lambda'\qquad\text{if $\lambda'_i-1\leq\lambda_i\leq\lambda'_i$ for each $i=1,\ldots,m$.}
\]
Let $\cals$ denote the set of equivalence classes of symbols.
We define several relations on $\cals$:
\begin{align*}
\calb^+ &=\{\,(\Lambda,\Lambda')\in\cals\times\cals
\mid\Upsilon(\Lambda_*)^\rmT\preccurlyeq\Upsilon(\Lambda'^*)^\rmT,\ \Upsilon(\Lambda'_*)^\rmT\preccurlyeq\Upsilon(\Lambda^*)^\rmT\,\};\\
\calb^- &=\{\,(\Lambda,\Lambda')\in\cals\times\cals
\mid\Upsilon(\Lambda^*)^\rmT\preccurlyeq\Upsilon(\Lambda'_*)^\rmT,\ \Upsilon(\Lambda'^*)^\rmT\preccurlyeq\Upsilon(\Lambda_*)^\rmT\,\};\\
\calb^+_{\rmU,\rmU} &=\left\{\,(\Lambda,\Lambda')\in\calb^+\mid
{\rm def}(\Lambda')=\begin{cases}
0, & \text{if ${\rm def}(\Lambda)=0$};\\
-{\rm def}(\Lambda)+1, & \text{if ${\rm def}(\Lambda)\neq 0$}
\end{cases}\,\right\};\\
\calb^-_{\rmU,\rmU} &=\left\{\,(\Lambda,\Lambda')\in\calb^-\mid
{\rm def}(\Lambda')=-{\rm def}(\Lambda)-1\,\right\};\\
\calb_{\rmU_n,\rmU_{n'}} &=
\begin{cases}
\{(\Lambda_\lambda,\Lambda_{\lambda'})\in\calb^+_{\rmU,\rmU}\mid |\lambda|=n,\ |\lambda'|=n'\,\}, & \text{if $n+n'$ is even};\\
\{(\Lambda_\lambda,\Lambda_{\lambda'})\in\calb^-_{\rmU,\rmU}\mid|\lambda|=n,\ |\lambda'|=n'\,\}, & \text{if $n+n'$ is odd}.
\end{cases}
\end{align*}

From above definitions and (\ref{0214})
we know that $(\Lambda_\lambda,\Lambda_{\lambda'})\in\calb_{\rmU_n,\rmU_{n'}}$ implies
either
\[
\ell(\lambda_\infty)=\ell(\lambda'_\infty)=0\quad\text{ or }\quad|\ell(\lambda_\infty)-\ell(\lambda'_\infty)|=1.
\]
The following proposition on the Howe correspondence of unipotent characters for a dual pair of
two unitary groups is rephrased from \cite{amr} th\'eor\`eme 5.15.
Note that here we need only to assume that the characteristic of the base field is not equal to $2$
(\cf.~\cite{pan-Lusztig-correspondence} proposition~5.13).

\begin{prop}\label{0307}
Let $(\bfG,\bfG')=(\rmU_n,\rmU_{n'})$ be a reductive dual pair of two unitary groups.
Then the decomposition of the unipotent part of the Weil character for the dual pair $(\bfG,\bfG')$
is given by
\[
\omega_{\bfG,\bfG',1}=\sum_{(\Lambda_\lambda,\Lambda_{\lambda'})\in\calb_{\bfG,\bfG'}}\rho_\lambda\otimes\rho_{\lambda'}.
\]
\end{prop}

For a finite classical group $G$,
let $\cale(G)$ (resp.~$\cale(G)_1$) denote the set of irreducible characters (resp.~unipotent characters) of $G$.
The proposition establishes a relation between $\cale(G)_1$ and $\cale(G')_1$ which will be called the
(unipotent part of the) \emph{$\Theta$-correspondence}.

For $\rho_\lambda\in\cale(G)_1$ or $\Lambda\in\cals_\bfG$,
we define
\begin{align*}
\Theta_{\bfG'}(\rho_\lambda)
&= \{\,\rho_{\lambda'}\in\cale(G')_1\mid(\Lambda_\lambda,\Lambda_{\lambda'})\in\calb_{\bfG,\bfG'}\,\}; \\
\Theta_{\bfG'}(\Lambda)
&= \{\,\Lambda'\in\cals_{\bfG'}\mid(\Lambda,\Lambda')\in\calb_{\bfG,\bfG'}\,\}.
\end{align*}
For $k\geq 0$, we define
\[
\Theta_{\bfG'}(\Lambda)_k =\begin{cases}
\{\,\Lambda'\in\Theta_{\bfG'}(\Lambda)\mid |\Upsilon(\Lambda')_*|=|\Upsilon(\Lambda)^*|-k\,\}, & \text{if $n+n'$ is even};\\
\{\,\Lambda'\in\Theta_{\bfG'}(\Lambda)\mid |\Upsilon(\Lambda')^*|=|\Upsilon(\Lambda)_*|-k\,\}, & \text{if $n+n'$ is odd}.
\end{cases}
\]
It is known that $\Theta_{\bfG'}(\Lambda)_k=\emptyset$ if $k$ is large enough and
\[
\Theta_{\bfG'}(\Lambda)=\bigsqcup_{k\geq 0}\Theta_{\bfG'}(\Lambda)_k.
\]

\subsection{Definition of $\theta_k$}\label{0308}
Consider the dual pair $(\bfG,\bfG')=(\rmU_n,\rmU_{n'})$.
Suppose that $\lambda\in\calp(n)$, $d=\ell(\lambda_\infty)$.
Then $\Lambda_\lambda\in\cals_{n,\delta}$ where $\delta=d$ if $d$ is even; $\delta=-d$ if $d$ is odd,
and $\Upsilon(\Lambda_\lambda)\in\calp_2(\frac{1}{2}(n-\frac{d(d+1)}{2}))$.
Similarly, if $\lambda'\in\calp(n')$, $d'=\ell(\lambda'_\infty)$,
then $\Upsilon(\Lambda_{\lambda'})\in\calp_2(\frac{1}{2}(n'-\frac{d'(d'+1)}{2}))$.
From Proposition~\ref{0307}, we have the following commutative diagram
\[
\begin{CD}
\cals_{n,\delta}  @> \Theta_{\bfG'} >> \cals_{n',\delta'} \\
@V \Upsilon VV @VV \Upsilon V \\
\calp_2(n-\tfrac{d(d+1)}{2}) @>>> \calp_2(n'-\tfrac{d'(d'+1)}{2})
\end{CD}
\]
where $\Upsilon$ is a bijection and $\Theta_{\bfG'}$ is a correspondence.
By abusing the notation, the correspondence in the bottom of the above diagram is also
denoted by $\Theta_{\bfG'}$.
Now we define
\begin{equation}\label{0301}
\tau=\frac{1}{2}\left[\left(n'-\frac{d'(d'+1)}{2}\right)-\left(n-\frac{d(d+1)}{2}\right)\right]
\end{equation}
and we want to construct a mapping
\[
\theta_k\colon\calp_2(\tfrac{1}{2}(n-\tfrac{d(d+1)}{2}))
\longrightarrow\calp_2(\tfrac{1}{2}(n'-\tfrac{d'(d'+1)}{2}))
\]
as follows when $\tau\geq 0$.
Note that the definition of $\theta_0$ is modified from \cite{akp} definition 5.
  
\begin{enumerate}
\item Suppose that $n+n'$ is even and $\tau\geq 0$.
Moreover, suppose also that
\begin{equation}\label{0602}
d'=
\begin{cases}
d, & \text{if $d=0$};\\
d-1, & \text{if $d$ is even and $d>0$};\\
d+1, & \text{if $d$ is odd}.
\end{cases}
\end{equation}
Let $\sqbinom{\mu_1,\ldots,\mu_{m_1}}{\nu_1,\ldots,\nu_{m_2}}\in\calp_2(\tfrac{1}{2}(n-\tfrac{d(d+1)}{2}))$.
For $0\leq k\leq\mu_1$, we define
\[
\theta_k\colon\sqbinom{\mu_1,\ldots,\mu_{m_1}}{\nu_1,\ldots,\nu_{m_2}}
\mapsto\sqbinom{\nu_1,\ldots,\nu_{m_2}}{\mu_2,\ldots,\mu_{m_1}}\cup\sqbinom{\tau+k}{\mu_1-k}.
\]

\item
Suppose that $n+n'$ is odd and $\tau\geq 0$.
Moreover, suppose that
\[
d'=
\begin{cases}
d+1, & \text{if $d$ is even;}\\
d-1, & \text{if $d$ is odd}.
\end{cases}
\]
Let $\sqbinom{\mu_1,\ldots,\mu_{m_1}}{\nu_1,\ldots,\nu_{m_2}}\in\calp_2(\tfrac{1}{2}(n-\tfrac{d(d+1)}{2}))$.
For $0\leq k\leq\nu_1$, we define
\[
\theta_k\colon\sqbinom{\mu_1,\ldots,\mu_{m_1}}{\nu_1,\ldots,\nu_{m_2}}
\mapsto\sqbinom{\nu_2,\ldots,\nu_{m_2}}{\mu_1,\ldots,\mu_{m_1}}\cup\sqbinom{\nu_1-k}{\tau+k}.
\]
\end{enumerate}

Via the bijection $\Upsilon$ in (\ref{0212}),
we will also regard $\theta_k$ as a mapping from
$\cals_{n,\delta}$ to $\cals_{n',\delta'}$ (when $\tau\geq 0$).

\begin{lem}\label{0304}
Let $\Lambda\in\cals_\bfG$ and $k\geq 0$.
Then $\theta_k(\Lambda)$ is the unique element of maximal order in $\Theta_{\bfG'}(\Lambda)_k$.
\end{lem}
\begin{proof}
The proof is similar to that of lemma~4.10 in \cite{pan-eta}.
\end{proof}

\begin{lem}\label{0303}
Let $\Lambda\in\cals_\bfG$.
There exists a unique index $k_0$ such that
\[
\ord(\theta_0(\Lambda))<\ord(\theta_1(\Lambda))<\cdots<\ord(\theta_{k_0-1}(\Lambda))
<\ord(\theta_{k_0}(\Lambda))>\ord(\theta_{k_0+1}(\Lambda))>\cdots,
\]
i.e., $\theta_{k_0}(\Lambda)$ is the unique element of maximal order in the set
$\{\,\theta_k(\Lambda)\mid k\geq 0\,\}$.
\end{lem}
\begin{proof}
Write $\Upsilon(\Lambda)=\sqbinom{\mu_1,\ldots,\mu_{m_1}}{\nu_1,\ldots,\nu_{m_2}}$.
First suppose that $n+n'$ is even.
If $\mu_1=0$, then $\{\,\theta_k(\Lambda)\mid k=0,\ldots,\mu_1\,\}=\{\theta_0(\Lambda)\}$
and there is nothing to prove.
So we may assume that $\mu_1\geq 1$.
By the similar argument in the proof of lemma~4.11 in \cite{pan-eta},
we know that from $\theta_k(\Lambda)$ to $\theta_{k+1}(\Lambda)$,
there is a unique entry $\alpha_k$ in the first row of $\theta_k(\Lambda)$ is
changed to $\alpha_k+2$ and all other entries in the first row are unchanged;
and there is a unique entry $\beta_k$ in the second row of $\theta_k(\Lambda)$ is
changed to $\beta_k-2$ and all other entries in the second row are unchanged.
Moreover, we know that
\begin{itemize}
\item all elements in the sequence $\langle\alpha_k\rangle$ are of the same parity
and the sequence is strictly increasing;

\item all elements in the sequence $\langle\beta_k\rangle$ are of the same parity
and the sequence is strictly decreasing;

\item any two elements $\alpha_k,\beta_{k'}$ are of opposite parties.
\end{itemize}
Now we define the index $k_0$ according to the following situations:
\begin{enumerate}
\item if $\alpha_0>\beta_0$, we let $k_0=0$;

\item if $\alpha_{\mu_1}<\beta_{\mu_1}$, we let $k_0=\mu_1$; or

\item if there is a unique index $k_1$ such that $\alpha_{k_1-1}<\beta_{k_1-1}$ and
$\alpha_{k_1}>\beta_{k_1}$, then we let
\[
k_0=\begin{cases}
k_1, & \text{if $\alpha_{k_1-1}+2<\beta_{k_1-1}$};\\
k_1-1, & \text{if $\alpha_{k_1-1}+2>\beta_{k_1-1}$}.
\end{cases}
\]
\end{enumerate}

Next we suppose that $n+n'$ is odd.
From above, we know that
from $\theta_k(\Lambda)$ to $\theta_{k+1}(\Lambda)$,
there is a unique entry $\alpha_k$ in the first row of $\theta_k(\Lambda)$ is
changed to $\alpha_k-2$ and all other entries in the first row are unchanged;
and there is a unique entry $\beta_k$ in the second row of $\theta_k(\Lambda)$ is
changed to $\beta_k+2$ and all other entries in the second row are unchanged.
Moreover,
\begin{itemize}
\item all elements in the sequence $\langle\alpha_k\rangle$ are of the same parity
and the sequence is strictly decreasing;

\item all elements in the sequence $\langle\beta_k\rangle$ are of the same parity
and the sequence is strictly increasing;

\item any two elements $\alpha_k,\beta_{k'}$ are of opposite parties.
\end{itemize}
Now we define the index $k_0$ according to the following situations:
\begin{enumerate}
\item if $\beta_0>\alpha_0$, we let $k_0=0$;

\item if $\beta_{\nu_1}<\alpha_{\nu_1}$, we let $k_0=\nu_1$; or

\item if there is a unique index $k_1$ such that $\beta_{k_1-1}<\alpha_{k_1-1}$ and
$\beta_{k_1}>\alpha_{k_1}$, then we let
\[
k_0=\begin{cases}
k_1, & \text{if $\beta_{k_1-1}+2<\alpha_{k_1-1}$};\\
k_1-1, & \text{if $\beta_{k_1-1}+2>\alpha_{k_1-1}$}.
\end{cases}
\]
\end{enumerate}

Then the assertion follows from Lemma~\ref{0210} immediately.
\end{proof}

\begin{exam}
Consider the dual pair $(\rmU_8,\rmU_{10})$.
Let $\lambda=[6,2]\in\calp(8)$.
Then $\lambda_\infty=[0]$, $X_\lambda=\{7,2\}$, $\Lambda=\Lambda_\lambda=\binom{7}{2}\in\cals_{8,0}$,
and $\Upsilon(\Lambda)=\sqbinom{3}{1}\in\calp_2(4)$.
Now the mappings $\theta_k\colon\calp_2(4)\rightarrow\calp_2(5)$ and $\theta_k\colon\cals_{8,0}\rightarrow\cals_{10,0}$
are given by
\begin{align*}
\textstyle\theta_0(\sqbinom{3}{1}) &=\textstyle\sqbinom{1,1}{3}, &
\textstyle\theta_1(\sqbinom{3}{1}) &=\textstyle\sqbinom{2,1}{2}, &
\textstyle\theta_2(\sqbinom{3}{1}) &=\textstyle\sqbinom{3,1}{1}, &
\textstyle\theta_3(\sqbinom{3}{1}) &=\textstyle\sqbinom{41,1}{0}; \\
\textstyle\theta_0(\binom{7}{2}) &=\textstyle\binom{5,3}{8,0}, &
\textstyle\theta_1(\binom{7}{2}) &=\textstyle\binom{7,3}{6,0}, &
\textstyle\theta_2(\binom{7}{2}) &=\textstyle\binom{9,3}{4,0}, &
\textstyle\theta_3(\binom{7}{2}) &=\textstyle\binom{11,3}{2,0}.
\end{align*}
Then the sequence $\langle\alpha_k\rangle$ is $5,7,9,11$, and
the sequence $\langle\beta_k\rangle$ is $8,6,4,2$.
Now $\alpha_0<\beta_0$, $\alpha_1>\beta_1$, and $\alpha_0+2<\beta_0$.
So we have $k_0=1$, i.e.,
\[
\textstyle
\ord(\theta_0(\binom{7}{2}))< \ord(\theta_1(\binom{7}{2}))> \ord(\theta_2(\binom{7}{2}))> \ord(\theta_3(\binom{7}{2})).
\]
\end{exam}

\subsection{Definitions of correspondences $\underline\theta$ and $\overline\theta$}\label{0309}
Keep the setting in the previous subsection, in particular, we assume that $\tau\geq 0$.
For $\Lambda\in\cals_{n,\delta}\subset\cals_\bfG$,
we define
\[
\underline\theta_{\bfG'}(\Lambda)=\theta_0(\Lambda)\quad\text{ and }\quad
\underline\theta_{\bfG'}(\rho_\Lambda)=\rho_{\theta_0(\Lambda)}.
\]
Then we define a relation between $\cale(G)_1$ and $\cale(G')_1$:
\[
\underline\theta_{\bfG,\bfG'}
=\{\,(\rho_\Lambda,\rho_{\Lambda'})\in\cale(G)_1\times\cale(G')_1\mid\Lambda'=\underline\theta_{\bfG'}(\Lambda)\,\}.
\]

\begin{lem}
Let $(\bfG,\bfG')=(\rmU_n,\rmU_{n'})$, and let $\lambda\in\calp(n)$.
Then $\underline\theta_{\bfG'}(\rho_\lambda)\in\Theta_{\bfG'}(\rho_\lambda)$.
\end{lem}
\begin{proof}
We know that $\calb_{\rmU_n,\rmU_{n'}}\subset\calb^+$ if $n+n'$ is even;
and $\calb_{\rmU_n,\rmU_{n'}}\subset\calb^-$ if $n+n'$ is odd.
Then the lemma follows from Lemma~\ref{0306} immediately.
\end{proof}

To define $\overline\theta_{\bfG'}$, we need to introduce a linear order ``$<$''on the set
$\cals_{n,\delta}\subset\cals_{\bfG}$ as follows:
\begin{enumerate}
\item Suppose that $n+n'$ is even.
Let $\Lambda,\Lambda'\in\cals_{n,\delta}$.
We define that $\Lambda<\Lambda'$ if either
\begin{itemize}
\item $|\Upsilon(\Lambda)^*|<|\Upsilon(\Lambda')^*|$; or

\item $|\Upsilon(\Lambda)^*|=|\Upsilon(\Lambda')^*|$ and $\Upsilon(\Lambda)^*<\Upsilon(\Lambda')^*$ in
lexicographic order; or

\item $\Upsilon(\Lambda)^*=\Upsilon(\Lambda')^*$ and $\Upsilon(\Lambda)_*<\Upsilon(\Lambda')_*$ in
lexicographic order.
\end{itemize}

\item Suppose that $n+n'$ is odd.
Let $\Lambda,\Lambda'\in\cals_{n,\delta}$.
We define that $\Lambda<\Lambda'$ if either
\begin{itemize}
\item $|\Upsilon(\Lambda)_*|<|\Upsilon(\Lambda')_*|$; or

\item $|\Upsilon(\Lambda)_*|=|\Upsilon(\Lambda')_*|$ and $\Upsilon(\Lambda)_*<\Upsilon(\Lambda')_*$ in
lexicographic order; or

\item $\Upsilon(\Lambda)_*=\Upsilon(\Lambda')_*$ and $\Upsilon(\Lambda)^*<\Upsilon(\Lambda')^*$ in
lexicographic order.
\end{itemize}
\end{enumerate}

Now we define $\overline\theta_{\bfG'}(\Lambda)$ inductively as follows.
Assume that $\overline\theta_{\bfG'}(\Lambda')$ is defined for all $\Lambda'<\Lambda$
and consider the set
\[
\Theta_{\bfG'}^\flat(\Lambda)
=\Theta_{\bfG'}(\Lambda)\smallsetminus\{\,\overline\theta_{\bfG'}(\Lambda')\mid\Lambda'<\Lambda\,\}.
\]

\begin{lem}
Let $\Lambda\in\cals_{n,\delta}\subset\cals_\bfG$.
Then the set $\Theta_{\bfG'}^\flat(\Lambda)$ is always non-empty.
\end{lem}
\begin{proof}
First suppose that $n+n'$ is even and $\Upsilon(\Lambda)=\sqbinom{\lambda_1,\ldots,\lambda_{m_1}}{\mu_1,\ldots,\mu_{m_2}}$.
Let $\Lambda'_0\in\Theta_{\bfG'}(\Lambda)_0$ be given such that
\[
\Upsilon(\Lambda'_0)=\sqbinom{\mu_1+\tau,\mu_2,\ldots,\mu_{m_2}}{\lambda_1,\ldots,\lambda_{m_1}}
\]
where $\tau$ is given in (\ref{0301}).
By the same argument in the proof of lemma~4.18 in \cite{pan-eta},
we see that $\Lambda'_0\not\in\Theta_{\bfG'}(\Lambda')$ for any $\Lambda'<\Lambda$,
and hence $\Lambda'_0$ is in $\Theta_{\bfG'}^\flat(\Lambda)$.

The proof for the case that $n+n'$ is odd is similar.
\end{proof}

Because now $\Theta^\flat_{\bfG'}(\Lambda)$ is non-empty,
we define $\overline\theta_{\bfG'}(\Lambda)$ to be the smallest element in the set of
elements of maximal order in $\Theta^\flat_{\bfG'}(\Lambda)$.
Then we have a mapping $\overline\theta_{\bfG'}\colon\cale(G)_1\rightarrow\cale(G')_1$ by
$\overline\theta_{\bfG'}(\rho_\Lambda)=\rho_{\overline\theta_{\bfG'}(\Lambda)}$,
and a relation $\overline\theta_{\bfG,\bfG'}$ between $\cale(G)_1$ and $\cale(G')_1$ by
\[
\overline\theta_{\bfG,\bfG'}
=\{\,(\rho_\Lambda,\rho_{\Lambda'})\in\cale(G)_1\times\cale(G')_1\mid\Lambda'=\overline\theta_{\bfG'}(\Lambda)\,\}.
\]

\begin{exam}
Consider the dual pair $(\rmU_7,\rmU_{10})$.
We know that
\[
\cals_{\rmU_7}=\cals_{7,-1}\cup\cals_{7,2}\quad\text{and}\quad
\cals_{\rmU_{10}}=\cals_{10,0}\cup\cals_{10,-3}\cup\cals_{10,4}.
\]
Now $\Upsilon$ establishes the bijections $\cals_{7,-1}\simeq\calp_2(3)$,
$\cals_{7,2}\simeq\calp_2(2)$, $\cals_{10,0}\simeq\calp_2(5)$,
$\cals_{10,-3}\simeq\calp_2(2)$ and $\cals_{10,4}\simeq\calp_2(0)$.
The correspondence $\Theta_{\bfG,\bfG'}$ between $\cals_{\rmU_7}$ and $\cals_{\rmU_{10}}$ is decomposed
as the union of the correspondence between $\cals_{7,-1}$ and $\cals_{10,0}$,
and the correspondence between $\cals_{7,2}$ and $\cals_{10,3}$.

Note that for the part of the correspondence $\cals_{7,-1}\rightarrow\cals_{10,0}$,
$\tau=5-3>0$,
so every element in $\cals_{7,-1}$ occurs in the correspondence $\Theta$.
Now we have the following table of the correspondence $\cals_{7,-1}\rightarrow\cals_{10,0}$.
A symbol $\Lambda'$ of maximal order in $\Theta_{\bfG'}(\Lambda)$ is superscripted by $\natural$
(Notation: ``$\Lambda'^\natural$''),
$\Lambda'$ is overlined (Notation: ``$\overline{\Lambda'}$'') if $\Lambda'=\overline\theta_{\bfG'}(\Lambda)$,
$\Lambda'\in\Theta_{\bfG'}(\Lambda)$ is cancelled out (Notation: ``$\bcancel{\Lambda'}$'') if
$\Lambda'\not\in\Theta^\flat_{\bfG'}(\Lambda)$.
The first element in $\Theta_{\bfG'}(\Lambda)_k$ is $\theta_k(\Lambda)$.
In particular, the first element in $\Theta_{\bfG'}(\Lambda)_0$ is $\theta_0(\Lambda)=\underline\theta_{\bfG'}(\Lambda)$.
\[
\begin{tabular}{c|llll}
\toprule
$\cals_{7,-1}$ & $\cals_{10,0}$ \\
$\Lambda$ & $\Theta_{\bfG'}(\Lambda)_0$ & $\Theta_{\bfG'}(\Lambda)_1$ & $\Theta_{\bfG'}(\Lambda)_2$ & $\Theta_{\bfG'}(\Lambda)_3$ \\
\midrule
$\binom{6,4,2}{7,5,3,1}$ & $\overline{\binom{7,5,3,1}{10,6,4,2}}^\natural,\binom{5,3,1}{10,4,2}$ & \\
$\binom{6,2}{5,3,1}$ & $\overline{\binom{5,3,1}{8,6,2}}^\natural,\binom{5,3,1}{10,4,2},\binom{3,1}{8,4},\binom{3,1}{10,2}$ & \\
$\binom{6}{3,1}$ & $\overline{\binom{3,1}{8,4}}^\natural,\binom{3,1}{10,2},\binom{1}{10}$ & \\
\midrule
$\binom{4,2}{7,3,1}$ & $\overline{\binom{7,3,1}{8,4,2}}^\natural,\binom{5,1}{8,2}$ & $\binom{5,3,1}{10,4,2},\binom{3,1}{10,2}$ & \\
$\binom{4}{5,1}$ & $\overline{\binom{5,1}{6,4}}^\natural,\binom{3}{8}$
& $\bcancel{\binom{3,1}{8,4}},\binom{3,1}{10,2},\binom{1}{10}$ & \\
\midrule
$\binom{2}{5,3}$ & $\overline{\binom{5,3}{6,2}}^\natural,\binom{5,3}{8,0}$ & $\binom{5,1}{8,2},\binom{3}{8}$ & \\
$\binom{2}{7,1}$ & $\overline{\binom{7,1}{6,2}}^\natural,\binom{5}{6}$ & $\binom{5,1}{8,2},\binom{3}{8}$ & $\binom{3,1}{10,2},\binom{1}{10}$ & \\
\midrule
$\binom{2,0}{7,5,3}$ & $\overline{\binom{7,5,3}{8,2,0}}^\natural$ & $\binom{5,3}{8,0}$ & \\
$\binom{0}{7,3}$ & $\overline{\binom{7,3}{6,0}}^\natural$ & $\binom{5,3}{8,0},\binom{5}{6}$ & $\binom{3}{8}$ & \\
$\binom{-}{7}$ & $\binom{7}{4}$ & $\overline{\binom{5}{6}}^\natural$ & $\binom{3}{8}$ & $\binom{1}{10}$ \\
\bottomrule
\end{tabular}
\]
Note that $\underline\theta_{\bfG'}(\binom{-}{7})=\binom{7}{4}\neq\binom{5}{6}=\overline\theta_{\bfG'}(\binom{-}{7})$,
and $\underline\theta_{\bfG'}(\Lambda)=\overline\theta_{\bfG'}(\Lambda)$ for
any other $\Lambda$ in $\cals_{7,-1}$.

Similarly, for the part $\cals_{7,2}\rightarrow\cals_{10,-3}$,
$\tau=2-2=0$,
so every element in $\cals_{7,2}$ and every elements in $\cals_{10,-3}$ occur in the correspondence $\Theta$.
Now we have the following table of the correspondence $\cals_{7,2}\rightarrow\cals_{10,-3}$:
\[
\begin{tabular}{c|llll}
\toprule
$\cals_{7,2}$ & $\cals_{10,-3}$ \\
$\Lambda$ & $\Theta_{\bfG'}(\Lambda)_0$ & $\Theta_{\bfG'}(\Lambda)_1$ & $\Theta_{\bfG'}(\Lambda)_2$ \\
\midrule
$\binom{5,3}{-}$ & $\overline{\binom{-}{7,5,1}}^\natural$ & \\
$\binom{7,1}{-}$ & $\overline{\binom{-}{9,3,1}}^\natural$ & \\
\midrule
$\binom{7,3,1}{2}$ & $\overline{\binom{2}{9,5,3,1}}^\natural$ & $\bcancel{\binom{-}{7,5,1}},\bcancel{\binom{-}{9,3,1}}$ & \\
\midrule
$\binom{7,5,3,1}{4,2}$ & $\overline{\binom{4,2}{9,7,5,3,1}}^\natural$ & $\bcancel{\binom{2}{9,5,3,1}}$ & \\
$\binom{5,3,1}{4}$ & $\overline{\binom{4}{7,5,3,1}}^\natural$ & $\bcancel{\binom{2}{9,5,3,1}}$ & $\bcancel{\binom{-}{9,3,1}}$ & \\
\bottomrule
\end{tabular}
\]
For this case we have $\underline\theta_{\bfG'}(\Lambda)=\overline\theta_{\bfG'}(\Lambda)$
for any $\Lambda\in\cals_{7,2}$.
This will be seen in Proposition~\ref{0302}.
\end{exam}

\begin{exam}
Consider the dual pair $(\rmU_8,\rmU_{10})$.
We know that
\[
\cals_{\rmU_8}=\cals_{8,0}\cup\cals_{8,-3}\quad\text{and}\quad
\cals_{\rmU_{10}}=\cals_{10,0}\cup\cals_{10,-3}\cup\cals_{10,4}.
\]
Now $\Upsilon$ establishes the bijections $\cals_{8,0}\simeq\calp_2(4)$,
$\cals_{8,-3}\simeq\calp_2(1)$, $\cals_{10,0}\simeq\calp_2(5)$,
$\cals_{10,-3}\simeq\calp_2(2)$ and $\cals_{10,4}\simeq\calp_2(0)$.
Now the correspondence $\Theta_{\bfG,\bfG'}$ between $\cals_{\rmU_8}$ and $\cals_{\rmU_{10}}$ is decomposed
as the union of the correspondence between $\cals_{8,0}$ and $\cals_{10,0}$ and the correspondence
between $\cals_{8,-3}$ and $\cals_{10,4}$.

Note that for the part of the correspondence $\cals_{8,0}\rightarrow\cals_{10,0}$,
$\tau=\frac{1}{2}(10-8)>0$, so every element in $\cals_{8,0}$ occurs in the correspondence $\Theta$.
Now we have the following table of the correspondence $\cals_{8,0}\rightarrow\cals_{10,0}$:
\[
\begin{tabular}{c|lllll}
\toprule
$\cals_{8,0}$ & $\cals_{10,0}$ \\
$\Lambda$ & $\Theta_{\bfG'}(\Lambda)_0$ & $\Theta_{\bfG'}(\Lambda)_1$ & $\Theta_{\bfG'}(\Lambda)_2$ & $\Theta_{\bfG'}(\Lambda)_3$ & $\Theta_{\bfG'}(\Lambda)_4$ \\
\midrule
$\binom{7,5,3,1}{8,6,4,2}$ & $\overline{\binom{11,9,7,5,3}{8,6,4,2,0}}^\natural,\binom{11,7,5,3}{6,4,2,0}$ & \\
$\binom{5,3,1}{8,4,2}$ & $\overline{\binom{11,7,5,3}{6,4,2,0}}^\natural,\binom{9,7,3}{4,2,0},\binom{11,5,3}{4,2,0}$ & \\
$\binom{3,1}{6,4}$ & $\overline{\binom{9,7,3}{4,2,0}}^\natural,\binom{9,5}{2,0}$ & \\
$\binom{3,1}{8,2}$ & $\overline{\binom{11,5,3}{4,2,0}}^\natural,\binom{9,5}{2,0},\binom{11,3}{2,0}$ & \\
$\binom{1}{8}$ & $\overline{\binom{11,3}{2,0}}^\natural,\binom{11}{0}$ & \\
\midrule
$\binom{7,3,1}{6,4,2}$ & $\overline{\binom{9,7,5,3}{8,4,2,0}}^\natural,\binom{9,5,3}{6,2,0}$ & $\bcancel{\binom{11,7,5,3}{6,4,2,0}},\bcancel{\binom{11,5,3}{4,2,0}}$ & \\
$\binom{5,1}{6,2}$ & $\overline{\binom{9,5,3}{6,2,0}}^\natural,\binom{7,5}{4,0},\binom{9,3}{4,0}$
& $\bcancel{\binom{9,7,3}{4,2,0}},\bcancel{\binom{11,5,3}{4,2,0}},\binom{9,5}{2,0},\bcancel{\binom{11,3}{2,0}}$ & \\
$\binom{3}{6}$ & $\overline{\binom{9,3}{4,0}}^\natural,\binom{9}{2}$ & $\binom{9,5}{2,0},\bcancel{\binom{11,3}{2,0}},\binom{11}{0}$ & \\
\midrule
$\binom{5,3}{4,2}$ & $\overline{\binom{7,5,3}{6,4,0}}^\natural,\binom{7,3}{4,2}$ & $\bcancel{\binom{9,5,3}{6,2,0}},\bcancel{\binom{9,3}{4,0}}$ & \\
$\binom{5,3}{6,0}$ & $\overline{\binom{7,3}{4,2}}^\natural,\binom{9,1}{4,2}$ & $\binom{7,5}{4,0},\bcancel{\binom{9,3}{4,0}},\binom{9}{2}$ & \\
$\binom{7,1}{4,2}$ & $\overline{\binom{7,5,3}{8,2,0}}^\natural,\binom{7,3}{6,0}$ & $\bcancel{\binom{9,4,2}{6,2,0}}^\natural,\bcancel{\binom{9,3}{4,0}}$
& $\bcancel{\binom{11,5,3}{4,2,0}},\bcancel{\binom{11,3}{2,0}}$ & \\
$\binom{5}{4}$ & $\binom{7,3}{6,0},\binom{7}{4}$ & $\overline{\binom{7,5}{4,0}}^\natural,\bcancel{\binom{9,3}{4,0}},\binom{9}{2}$
& $\binom{9,5}{2,0},\bcancel{\binom{11,3}{2,0}},\binom{11}{0}$ & \\
\midrule
$\binom{7,5,3}{6,2,0}$ & $\overline{\binom{7,5,1}{6,4,2}}^\natural,\binom{9,3,1}{6,4,2}$ & $\bcancel{\binom{7,3}{4,2}},\binom{9,1}{4,2}$ & \\
$\binom{7,3}{4,0}$ & $\overline{\binom{5,3}{6,2}}^\natural,\binom{7,1}{6,2}$ & $\bcancel{\binom{7,3}{4,2}},\binom{9,1}{4,2},\binom{7,3}{6,0},\binom{7}{4}$
& $\bcancel{\binom{9,3}{4,0}},\binom{9}{2}$ & \\
$\binom{7}{2}$ & $\binom{5,3}{8,0},\binom{5}{6}$ & $\overline{\binom{7,3}{6,0}}^\natural,\binom{7}{4}$ & $\bcancel{\binom{9,3}{4,0}},\binom{9}{2}$
& $\bcancel{\binom{11,3}{2,0}},\binom{11}{0}$ & \\
\midrule
$\binom{9,7,5,3}{6,4,2,0}$ & $\overline{\binom{9,5,3,1}{8,6,4,2}}^\natural$ & $\binom{9,3,1}{6,4,2}$ & \\
$\binom{9,5,3}{4,2,0}$ & $\overline{\binom{7,3,1}{8,4,2}}^\natural$ & $\binom{9,3,1}{6,4,2},\binom{7,1}{6,2}$ & $\binom{9,1}{4,2}$ & \\
$\binom{7,5}{2,0}$ & $\overline{\binom{5,1}{6,4}}^\natural$ & $\binom{7,1}{6,2}$ & $\binom{7}{4}$ & \\
$\binom{9,3}{2,0}$ & $\overline{\binom{5,1}{8,2}}^\natural$ & $\binom{7,1}{6,2}^\natural,\binom{5}{6}$ & $\binom{9,1}{4,2}$ & $\binom{9}{2}$ & \\
$\binom{9}{0}$ & $\binom{3}{8}$ & $\overline{\binom{5}{6}}^\natural$ & $\binom{7}{4}$ & $\binom{9}{2}$ & $\binom{11}{0}$ \\
\bottomrule
\end{tabular}
\]
From the above table we see that for most $\Lambda\in\cals_{8,0}$,
we have $\underline\theta_{\bfG'}(\Lambda)=\overline\theta_{\bfG'}(\Lambda)$.
However, $\underline\theta_{\bfG'}(\Lambda)=\theta_0(\Lambda)\neq\theta_1(\Lambda)=\overline\theta_{\bfG'}(\Lambda)$
for $\Lambda=\binom{5}{4},\binom{7}{2},\binom{9}{0}$.

For the part $\cals_{8,-3}\rightarrow\cals_{10,4}$,
we have $\tau=1-2<0$,
so not every element in $\cals_{8,-3}$ occurs in the correspondence $\Theta$.
Then we have to switch the roles of $\bfG$ and $\bfG'$ and consider the correspondence $\cals_{10,4}\rightarrow\cals_{8,-3}$.
Note that $\binom{7,5,3,1}{-}$ is the only element in $\cals_{10,4}$,
and $\binom{-}{7,3,1},\binom{2}{7,5,3,1}$ are the two elements in $\cals_{8,-3}$.
The following table is the correspondence $\cals_{10,4}\rightarrow\cals_{8,-3}$:
\[
\begin{tabular}{c|lllll}
\toprule
$\cals_{10,4}$ & $\cals_{8,3}$ \\
$\Lambda'$ & $\Theta_\bfG(\Lambda')_0$ \\
\midrule
$\binom{7,5,3,1}{-}$ & $\overline{\binom{2}{7,5,3,1}}^\natural$ \\
\bottomrule
\end{tabular}
\]
Note that the symbol $\binom{-}{7,3,1}$ is the only element in $\cals_{\rmU_8}$ which does not occur
in the correspondence $\Theta_{\bfG,\bfG'}$.
\end{exam}

\subsection{Properties of correspondences $\underline\theta$ and $\overline\theta$}
Let $(\bfG,\bfG')=(\rmU_n,\rmU_{n'})$, $\lambda\in\calp(n)$, and $d=\ell(\lambda_\infty)$.
\begin{enumerate}
\item Suppose that $n+n'$ is even and $d$ is even.
Then from (\ref{0602}), we know that $\rho_\lambda$ occurs in the $\underline\theta_{\bfG,\bfG'}$ if
\[
n'-\frac{d(d-1)}{2}\geq n-\frac{d(d+1)}{2},
\]
i.e., if $n'\geq n-d$.

\item Suppose that $n+n'$ is even and $d$ is odd.
Then from (\ref{0602}), we know that $\rho_\lambda$ occurs in the $\underline\theta_{\bfG,\bfG'}$ if
\[
n'-\frac{(d+1)(d+2)}{2}\geq n-\frac{d(d+1)}{2},
\]
i.e., if $n'\geq n+d+1$.

\item Suppose that $n+n'$ is odd and $d$ is even.
Then $\rho_\lambda$ occurs in the $\underline\theta_{\bfG,\bfG'}$ if $n'\geq n+d+1$.

\item Suppose that $n+n'$ is odd and $d$ is odd.
Then $\rho_\lambda$ occurs in the $\underline\theta_{\bfG,\bfG'}$ if $n'\geq n-d$.
\end{enumerate}

\begin{prop}\label{0302}
Let $\Lambda\in\cals_{n,\delta}\subset\cals_\bfG$ and suppose that $\tau=0$.
Then $\underline\theta_{\bfG'}(\Lambda)=\overline\theta_{\bfG'}(\Lambda)$.
\end{prop}
\begin{proof}
By the same argument in the proof of lemma~5.2 in \cite{pan-eta},
we can show that $\Theta^\flat_{\bfG'}(\Lambda)=\{\theta_0(\Lambda)\}$ when $\tau=0$.
This implies that $\overline\theta_{\bfG'}(\Lambda)=\theta_0(\Lambda)=\underline\theta_{\bfG'}(\Lambda)$.
\end{proof}

A dual pair $(\rmU_n,\rmU_{n'})$ is called in \emph{stable range} if $n\leq\lfloor\frac{n'}{2}\rfloor$.

\begin{lem}
Suppose that $(\bfG,\bfG')$ is in stable range and let $\Lambda\in\cals_\bfG$.
Then $\theta_0(\Lambda)$ is the unique element of maximal order in $\Theta_{\bfG'}(\Lambda)$.
\end{lem}
\begin{proof}
Suppose that $(\bfG,\bfG')=(\rmU_n,\rmU_{n'})$.
Write
\[
\Lambda=\binom{a_1,\ldots,a_{m_1}}{b_1,\ldots,b_{m_2}}\in\cals_{n,\delta}\subset\cals_\bfG\quad\text{and}\quad
\Upsilon(\Lambda)=\sqbinom{\mu_1,\ldots,\mu_{m_1}}{\nu_1,\ldots,\nu_{m_2}}
\in\calp_2(\tfrac{1}{2}(n-\tfrac{d(d+1)}{2}))
\]
where $d=|\delta|$.
Now we consider the following cases:
\begin{enumerate}
\item Suppose that both $n,n'$ are even and $n'\geq 2n$.
Note that now
\begin{align*}
\tau-\nu_1
&\geq \tfrac{1}{2}\left(n'-n+\tfrac{d(d+1)}{2}-\tfrac{d'(d'+1)}{2}\right)-\tfrac{1}{2}\left(n-\tfrac{d(d+1)}{2}\right) \\
&=\tfrac{1}{2}(n'-2n)+\left(d(d+1)-\tfrac{d'(d'+1)}{2}\right).
\end{align*}

\begin{enumerate}
\item If $\delta$ even and positive,
then $d=\delta=m_1-m_2\geq 4$, $d'=d-1$, and
\[
\tau-\nu_1\geq\tfrac{1}{2}(n'-2n)+\tfrac{1}{2}d(d+3)\geq 0.
\]
So now
\begin{equation}\label{0305}
\theta_0(\Upsilon(\Lambda))=\sqbinom{\tau,\nu_1,\ldots,\nu_{m_2}}{\mu_1,\ldots,\mu_{m_1}},
\end{equation}
$n'$ is even, and $d'$ is odd.
Hence
\begin{align*}
\alpha_0 &= 2(\tau+m_2)=n'-n+d+2m_2, \\
\beta_0 &\leq a_1=\mu_1+2m_1-1\leq n-\tfrac{d(d+1)}{2}+2m_1-1, \\
\alpha_0-\beta_0 &\geq n'-2n-d+\tfrac{d(d+1)}{2}+1=n'-2n+\tfrac{d(d-1)}{2}+1>0.
\end{align*}

\item If $\delta=0$,
then $d'=d=m_1-m_2=0$, and $\tau-\nu_1=\tfrac{1}{2}(n'-2n)\geq 0$.
So now $\theta_0(\Upsilon(\Lambda))$ is as in (\ref{0305}), and both $n',d'$ are even.
Hence
\begin{align*}
\alpha_0 &= 2(\tau+m_2)+1=n'-n+2m_2+1, \\
\beta_0 &\leq a_1\leq n+2(m_1-1), \\
\alpha_0-\beta_0 &\geq n'-2n+3>0.
\end{align*}

\item If $\delta$ odd,
then $d=m_2-m_1=-\delta\geq 3$, $d'=d+1$, and
\[
\tau-\nu_1\geq \tfrac{1}{2}(n'-2n)+\tfrac{1}{2}(d+1)(d-2)\geq 0.
\]
So now $\theta_0(\Upsilon(\Lambda))$ is as in (\ref{0305}), and both $n',d'$ are even.
Hence
\begin{align*}
\alpha_0 &= 2(\tau+m_2)+1=n'-n-(d+1)+2m_2+1, \\
\beta_0 &\leq a_1=n-\tfrac{d(d+1)}{2}+2(m_1-1), \\
\alpha_0-\beta_0 &\geq n'-2n+\tfrac{d(d+3)}{2}>0.
\end{align*}
\end{enumerate}

\item Suppose that both $n,n'$ are odd and $n'\geq 2n+1$.
\begin{enumerate}
\item If $\delta$ is even,
then $d=\delta=m_1-m_2\geq 2$, $d'=d-1$, and
\[
\tau-\nu_1\geq \tfrac{1}{2}(n'-2n)+\tfrac{1}{2}d(d+3)\geq 0.
\]
So now $\theta_0(\Upsilon(\Lambda))$ is as in (\ref{0305}), and both $n',d'$ are odd.
Hence
\begin{align*}
\alpha_0 &\geq 2(\tau+m_2-1)=n'-n+d+2m_2-2, \\
\beta_0 &\leq a_1=n-\tfrac{d(d+1)}{2}+2m_1-1, \\
\alpha_0-\beta_0 &\geq n'-2n-d+\tfrac{d(d+1)}{2}-1=n'-2n+\tfrac{d(d-1)}{2}-1>0.
\end{align*}

\item If $\delta$ is odd,
then $d=m_2-m_1=-\delta\geq 1$, $d'=d+1$, and
\[
\tau-\nu_1\geq \tfrac{1}{2}(n'-2n)+\tfrac{1}{2}(d+1)(d-2)\geq 0.
\]
So now $\theta_0(\Upsilon(\Lambda))$ is as in (\ref{0305}), $n'$ is odd, and $d'$ is even.
Hence
\begin{align*}
\alpha_0 &\geq 2(\tau+m_2-1)+1=n'-n-(d+1)+2m_2-1, \\
\beta_0 &\leq a_1=n-\tfrac{d(d+1)}{2}+2(m_1-1), \\
\alpha_0-\beta_0 &\geq n'-2n+\tfrac{d(d+3)}{2}>0.
\end{align*}
\end{enumerate}

\item Suppose that $n$ is even, $n'$ is odd and $n'\geq 2n+1$.
Note that now
\begin{align*}
\tau-\mu_1
&\geq \tfrac{1}{2}\left(n'-n+\tfrac{d(d+1)}{2}-\tfrac{d'(d'+1)}{2}\right)-\tfrac{1}{2}\left(n-\tfrac{d(d+1)}{2}\right) \\
&=\tfrac{1}{2}(n'-2n)+\left(d(d+1)-\tfrac{d'(d'+1)}{2}\right).
\end{align*}
\begin{enumerate}
\item If $\delta$ is even,
then $d=\delta=m_1-m_2\geq 0$, $d'=d+1$, and
\[
\tau-\mu_1\geq \tfrac{1}{2}(n'-2n)+\tfrac{1}{2}(d+1)(d-2)\geq 0.
\]
So now
\begin{equation}\label{0306}
\theta_0(\Upsilon(\Lambda))=\sqbinom{\nu_1,\ldots,\nu_{m_2}}{\tau,\mu_1,\ldots,\mu_{m_1}},
\end{equation}
both $n'$ and $d'$ are odd.
Hence
\begin{align*}
\beta_0 &= 2(\tau+m_1)+1=n'-n-(d+1)+2m_1+1, \\
\alpha_0 &\leq a_1\leq n-\tfrac{d(d+1)}{2}+2(m_2-1), \\
\beta_0-\alpha_0 &\geq n'-2n+d+2+\tfrac{d(d+1)}{2}=n'-2n+\tfrac{d(d+3)}{2}+2>0.
\end{align*}

\item If $\delta$ is odd,
then $d=m_2-m_1=-\delta\geq 3$, $d'=d-1$, and $\tau-\mu_1\geq 0$.
So now $\theta_0(\Upsilon(\Lambda))$ is as in (\ref{0306}), $n'$ is odd, and $d'$ is even.
Hence
\begin{align*}
\beta_0 &= 2(\tau+m_1)=n'-n+d+2m_1, \\
\alpha_0 &\leq a_1=n-\tfrac{d(d+1)}{2}+2(m_2-1)+1, \\
\beta_0-\alpha_0 &\geq n'-2n+\tfrac{d(d-1)}{2}+1>0.
\end{align*}
\end{enumerate}

\item Suppose that $n$ is odd, $n'$ is even and $n'\geq 2n$.
\begin{enumerate}
\item If $\delta$ is even,
then $d=\delta=m_1-m_2\geq 2$, $d'=d+1$, and $\tau-\mu_1\geq 0$.
So now $\theta_0(\Upsilon(\Lambda))$ is as in (\ref{0306}), $n'$ is even, and $d'$ is odd.
Hence
\begin{align*}
\beta_0 &\geq 2(\tau+m_1)+1=n'-n-d+2m_1, \\
\alpha_0 &\leq a_1=n-\tfrac{d(d+1)}{2}+2(m_2-1), \\
\beta_0-\alpha_0 &\geq n'-2n+d+\tfrac{d(d+1)}{2}+2=n'-2n+\tfrac{d(d+3)}{2}+2>0.
\end{align*}

\item If $\delta$ is odd,
then $d=m_2-m_1=-\delta\geq 1$, $d'=d-1$, and $\tau-\mu_1\geq 0$.
So now $\theta_0(\Upsilon(\Lambda))$ is as in (\ref{0306}), and both $n',d'$ are even.
Hence
\begin{align*}
\beta_0 &\geq 2(\tau+m_1)=n'-n+d+2m_1, \\
\alpha_0 &\leq a_1=n-\tfrac{d(d+1)}{2}+2m_2-1, \\
\beta_0-\alpha_0 &\geq n'-2n-d+\tfrac{d(d+1)}{2}+1=n'-2n+\tfrac{d(d-1)}{2}+1>0.
\end{align*}
\end{enumerate}
\end{enumerate}
So we conclude that $\alpha_0>\beta_0$ when $n+n'$ is even;
and $\beta_0>\alpha_0$ when $n+n'$ is odd.
Then from the proof of Lemma~\ref{0303},
we see that $\theta_0(\Lambda)$ is the unique element of maximal order in the set
$\{\,\theta_k(\Lambda)\mid k\geq 0\,\}$.
Then the lemma follows from Lemma~\ref{0304} immediately.
\end{proof}

\begin{prop}\label{0310}
Suppose that the dual pair $(\bfG,\bfG')=(\rmU_n,\rmU_{n'})$ is in stable range
and let $\Lambda\in\cals_\bfG$.
Then $\underline\theta_{\bfG'}(\Lambda)=\overline\theta_{\bfG'}(\Lambda)$.
\end{prop}
\begin{proof}
Suppose that $(\bfG,\bfG')$ is in stable range and $\Lambda\in\cals_{n,\delta}\subset\cale_\bfG$ for some $\delta$.
Note that the mapping $\theta_0\colon\cals_{n,\delta}\rightarrow\cals_{n',\delta'}$ where $\delta'$ is given
in Subsection~\ref{0308} is one-to-one.
This implies that $\theta_0(\Lambda)\in\Theta^\flat_{\bfG'}(\Lambda)$.
Then by the previous lemma,
we have $\overline\theta_{\bfG'}(\Lambda)=\theta_0(\Lambda)=\underline\theta_{\bfG'}(\Lambda)$.
\end{proof}

Recall that we define the mappings $\underline\theta,\overline\theta\colon\cals_{n,\delta}\rightarrow\cals_{n',\delta'}$
under the assumption $\tau\geq 0$ in Subsection~\ref{0309}.
Now we extend the domain of both mappings by defining
\begin{align*}
\underline\theta_{\bfG'}(\rho_\lambda)
&=\rho_{\lambda'}\quad\text{ if and only if }\quad\underline\theta_\bfG(\rho_{\lambda'})=\rho_\lambda; \\
\overline\theta_{\bfG'}(\rho_\lambda)
&=\rho_{\lambda'}\quad\text{ if and only if }\quad\overline\theta_\bfG(\rho_{\lambda'})=\rho_\lambda.
\end{align*}
So from now on, we will drop the assumption that $\tau\geq0$.
However, $\underline\theta_{\bfG'}(\rho_\lambda)$ or $\overline\theta_{\bfG'}(\rho_\lambda)$
might not be defined if $\tau<0$. 

\section{Maximal Theta Relations for Unitary Groups}

Let $(\bfG,\bfG')=(\rmU_n,\rmU_{n'})$ for some $n,n'\in\bbN\cup\{0\}$.

\subsection{Lusztig correspondence and $\Theta$-correspondence}
For a semisimple element $s$ in the dual group $G^*$ of $G$,
let $\cale(G)_s$ denote the \emph{Lusztig series} associated to $s$.
For $\rho\in\cale(G)_s$,
let $G^{(1)}$, $G^{(2)}$ be defined as in \cite{pan-Lusztig-correspondence} so that
$C_{G^*}(s)=G^{(1)}\times G^{(2)}$.
Then we have a bijection called a \emph{Lusztig correspondence}
\[
\Xi_s\colon \cale(G)_s\longrightarrow\cale(G^{(1)}\times G^{(2)})_1.
\]
Write $\Xi_s(\rho)=\rho^{(1)}\otimes\rho^{(2)}$ for $\rho^{(j)}\in\cale(G^{(j)})_1$ and $j=1,2$.
The following can be extracted from \cite{amr} th\'eor\`eme 2.6 (\cf.~\cite{pan-chain01} theorem 3.10).
Note that from \cite{pan-finite-unipotent}, we do not need to assume that $q$ is large enough.

\begin{prop}
Let $(G,G')=(\rmU_n(q),\rmU_{n'}(q))$.
Let $\eta\in\cale(G)_s$ and $\eta\in\cale(G')_{s'}$ for some $s,s'$.
Then $\rho\otimes\rho'$ occurs in the Howe correspondence for $(G,G')$ if and only if
the following conditions hold:
\begin{itemize}
\item $G^{(1)}\simeq G'^{(1)}$ and $\rho^{(1)}\simeq\rho'^{(1)}$,

\item $\rho^{(2)}\otimes\rho'^{(2)}$ occurs in the correspondence for the dual pair $(G^{(2)},G'^{(2)})$,
\end{itemize}
i.e., the following diagram
\[
\begin{CD}
\rho @> \Theta_{\bfG,\bfG'} >> \rho' \\
@V \Xi_s VV @VV \Xi_{s'} V \\
\rho^{(1)}\otimes\rho^{(2)} @> {\rm id}\otimes\Theta_{\bfG^{(2)},\bfG'^{(2)}} >> \rho'^{(1)}\otimes\rho'^{(2)}.
\end{CD}
\]
commutes.
\end{prop}

Then we define
\begin{align}\label{0401}
\begin{split}
\underline\theta_{\bfG'}(\rho)
&= \Xi_{s'}^{-1}(\rho^{(1)}\otimes\underline\theta_{\bfG^{(2)}}(\rho^{(2)})) \\
\overline\theta_{\bfG'}(\rho)
&= \Xi_{s'}^{-1}(\rho^{(1)}\otimes\overline\theta_{\bfG^{(2)}}(\rho^{(2)})),
\end{split}
\end{align}
i.e., we have the following two commutative diagrams:
\[
\begin{CD}
\rho @> \underline\theta_{\bfG,\bfG'} >> \rho' \\
@V \Xi_s VV @VV \Xi_{s'} V \\
\rho^{(1)}\otimes\rho^{(2)} @> {\rm id}\otimes\underline\theta_{\bfG^{(2)},\bfG'^{(2)}} >> \rho'^{(1)}\otimes\rho'^{(2)},
\end{CD}
\qquad\qquad
\begin{CD}
\rho @> \overline\theta_{\bfG,\bfG'} >> \rho' \\
@V \Xi_s VV @VV \Xi_{s'} V \\
\rho^{(1)}\otimes\rho^{(2)} @> {\rm id}\otimes\overline\theta_{\bfG^{(2)},\bfG'^{(2)}} >> \rho'^{(1)}\otimes\rho'^{(2)}.
\end{CD}
\]
Therefore the domains of the relations $\underline\theta$ and $\overline\theta$ are extended from unipotent
characters to all irreducible characters.

\begin{cor}
Suppose that the dual pair $(\bfG,\bfG')=(\rmU_n,\rmU_{n'})$ is in stable range.
Then $\underline\theta=\overline\theta$.
\end{cor}
\begin{proof}
This follows from Proposition~\ref{0310} and (\ref{0401}) immediately.
\end{proof}

For a fixed Witt tower,
let $\bfG_n$ denote the group of split rank $n$,
i.e., $\bfG_n=\rmU_{2n}$ or $\rmU_{2n+1}$.
For $\rho\in\cale(G)$,
let $n'_0(\rho)$ (resp.~$\underline n'_0(\rho)$) be the smallest $n'$ such that
$\Theta_{\bfG'_{n'}}(\rho)\neq\emptyset$ (resp.~$\underline\theta_{\bfG'_{n'}}(\rho)$ is defined).

\begin{cor}\label{0603}
Let $(\bfG,\bfG')$ be a dual pair of two unitary groups.
Suppose that $\rho'\in\cale(G')$.
\[
n_0(\rho')=\underline n_0(\rho').
\]
\end{cor}
\begin{proof}
The proof is similar to that of lemma~6.9 in \cite{pan-eta}.
\end{proof}

\subsection{Maximal theta relation}\label{0402}
Let $\vartheta$ be a sub-relation of $\Theta$ (\cf.~subsection 7.1 in \cite{pan-eta}).
For a dual pair $(\bfG,\bfG')$ and $\rho\in\cale(G)$,
we denote
\[
\vartheta_{\bfG'}(\rho)=\{\,\rho'\in\cale(G')\mid(\rho,\rho')\in\vartheta_{\bfG,\bfG'}\,\}.
\]
A partial ordering is given on the set of all sub-relations of $\Theta$ by inclusion, i.e.,
we say that $\vartheta_1\leq\vartheta_2$ if, for each dual pair $(\bfG,\bfG')$ and $\rho\in\cale(G)$,
we have $\vartheta_{1,\bfG'}(\rho)\subseteq\vartheta_{2,\bfG'}(\rho)$.

\begin{itemize}
\item A sub-relation $\vartheta$ of $\Theta$ is called \emph{semi-persistent} (on unipotent characters) if and $\rho_\lambda\in\cale(\rmU_n(q))_1$ occurs in $\vartheta_{\rmU_n,\rmU_{n'}}$ whenever either
\begin{itemize}
\item $n+n'+\ell(\lambda_\infty)$ is even, and $n'\geq n-\ell(\lambda_\infty)$; or

\item $n+n'+\ell(\lambda_\infty)$ is odd and $n'\geq n+\ell(\lambda_\infty)+1$.
\end{itemize}

\item A sub-relation $\vartheta$ is called \emph{symmetric} if for each dual pair
$(\bfG,\bfG')=(\rmU_n,\rmU_{n'})$, $\rho\in\cale(G)$ and $\rho'\in\cale(G')$,
we have $\rho'\in\vartheta_{\bfG'}(\rho)$ if and only if $\rho\in\vartheta_\bfG(\rho')$.

\item A sub-relation $\vartheta$ is said to be \emph{compatible with} the Lusztig correspondence if
for each dual pair $(\bfG,\bfG')=(\rmU_n,\rmU_{n'})$, $\rho\in\cale(G)$ and $\rho'\in\cale(G')$,
the following diagram
\[
\begin{CD}
\rho @> \vartheta_{\bfG,\bfG'} >> \rho' \\
@V \Xi_s VV @VV \Xi_{s'} V \\
\rho^{(1)}\otimes\rho^{(2)} @> {\rm id}\otimes\vartheta_{\bfG^{(2)},\bfG'^{(2)}} >> \rho'^{(1)}\otimes\rho'^{(2)}.
\end{CD}
\]
commutes.
\end{itemize}

Similar to the case for symplectic/orthogonal dual pair,
a sub-relation $\vartheta$ of $\Theta$ is called a \emph{theta relation} if it is
semi-persistent, symmetric and compatible with the Lusztig correspondence.
We know that both $\underline\theta$ and $\overline\theta$ are one-to-one theta-relations.
The following proposition says that a one-to-one theta-relation can not be properly contained
in another one-to-one theta-relation.

\begin{prop}
In the set of all one-to-one theta-relations,
each element is maximal.
\end{prop}
\begin{proof}
Let $\vartheta$ be a one-to-one theta-relation.
Suppose that $\vartheta'$ is another theta-relation such that $\vartheta'_{\bfG,\bfG'}$ properly contains $\vartheta_{\bfG,\bfG'}$ for some dual pair $(\bfG,\bfG')=(\rmU_n,\rmU_{n'})$, i.e.,
there are $\rho\in\cale(G)$ and $\rho'\in\cale(G')$ such that
$(\rho,\rho')\in\vartheta'_{\bfG,\bfG'}$ and $(\rho,\rho')\not\in\vartheta_{\bfG,\bfG'}$.
If $\vartheta_{\bfG'}(\rho)$ is defined,
then $\vartheta'$ is not one-to-one.
So we may assume that $\vartheta_{\bfG'}(\rho)$ is not defined.
Since both $\vartheta'$ and $\vartheta$ are compatible with Lusztig correspondence,
we may assume that both $\rho,\rho'$ are unipotent.
So we write $\rho=\rho_\lambda$ and $\rho'=\rho_{\lambda'}$ for some
$\lambda\in\calp(n)$ and $\lambda'\in\calp(n')$.

First suppose that $n+n'$ is even.
Because $\vartheta_{\bfG'}(\rho_\lambda)$ is not defined,
we must have
\[
n'<
\begin{cases}
n-\ell(\lambda_\infty), & \text{if $\ell(\lambda_\infty)$ is even};\\
n+\ell(\lambda_\infty)+1, & \text{if $\ell(\lambda_\infty)$ is odd}.
\end{cases}
\]
Because now $(\rho_\lambda,\rho_{\lambda'})\in\vartheta'_{\bfG,\bfG'}\subseteq\Theta_{\bfG,\bfG'}$,
we have $(\Lambda_\lambda,\Lambda_{\lambda})\in\calb_{\bfG,\bfG'}$.
Therefore, we have
\[
\ell(\lambda'_\infty)=
\begin{cases}
\ell(\lambda_\infty), & \text{if $\ell(\lambda_\infty)=0$};\\
\ell(\lambda_\infty)-1, & \text{if $\ell(\lambda_\infty)$ is even and $\ell(\lambda_\infty)>0$};\\
\ell(\lambda_\infty)+1, & \text{if $\ell(\lambda_\infty)$ is odd},
\end{cases}
\]
and then
\[
\ell(\lambda_\infty)=
\begin{cases}
\ell(\lambda'_\infty), & \text{if $\ell(\lambda'_\infty)=0$};\\
\ell(\lambda'_\infty)-1, & \text{if $\ell(\lambda'_\infty)$ is odd};\\
\ell(\lambda'_\infty)+1, & \text{if $\ell(\lambda'_\infty)$ is even and $\ell(\lambda_\infty)>0$}.
\end{cases}
\]
Hence we have
\[
n\geq
\begin{cases}
n'-\ell(\lambda'_\infty), & \text{if $\ell(\lambda'_\infty)$ is even};\\
n'+\ell(\lambda'_\infty)+1, & \text{if $\ell(\lambda'_\infty)$ is odd}.
\end{cases}
\]
So we have $\vartheta_{\bfG}(\rho_{\lambda'})$ is defined.
Let $\rho_{\lambda''}=\vartheta_{\bfG}(\rho_\lambda)\in\cale(G)$.
Now $(\rho_{\lambda'},\rho_{\lambda''})\in\vartheta_{\bfG',\bfG}\subseteq\vartheta'_{\bfG',\bfG}$ implies
that $(\rho_{\lambda''},\rho_{\lambda'})\in\vartheta'_{\bfG,\bfG'}$ since $\vartheta$ is symmetric.
Moreover, $(\rho_\lambda,\rho_{\lambda'})$ is in $\vartheta'_{\bfG,\bfG'}$ by our assumption.
However, $\lambda$ and $\lambda''$ are not equal
because $\vartheta_{\bfG'}(\rho_\lambda)$ is not defined and
$\vartheta_{\bfG'}(\rho_{\lambda''})=\rho_{\lambda'}$
by the symmetricity of $\underline\theta$.
So we conclude that $\vartheta'$ is not one-to-one.

The proof for the case that $n+n'$ is odd is similar.
\end{proof}

\begin{cor}
Let $(\bfG,\bfG')$ be a dual pair of two unitary groups.
Then both $\underline\theta$ and $\overline\theta$ are maximal one-to-one theta-relations.
\end{cor}

\bibliography{refer}
\bibliographystyle{amsalpha}

\end{document}